\documentclass[12pt]{amsart}
\usepackage{ifthen,color,xypic,amssymb}

\newcommand{\catD}{{D}}
\newcommand{\dbc}[1]{{D}^b(#1)}
\newcommand{\dpc}[1]{{D}^+(#1)}
\newcommand{\dmc}[1]{{D}^-(#1)}

\newcommand{\cdbc}[1]{{D}^b_c(#1)}
\newcommand{\fmf}[3]{{\Phi^{#1}_{{\scriptscriptstyle #2\!\rightarrow\! #3}}}}
\newcommand{\sk}[1]{{\cO_{#1}}}
\newcommand{\Hom}{{\operatorname{Hom}}}
\newcommand{\Ext}{{\operatorname{Ext}}}

\newcommand{\Tor}{{\operatorname{Tor}}}
\newcommand{\SHom}{{\mathcal{H}om}}

\newcommand{\lotimes}{{\,\stackrel{\mathbf L}{\otimes}\,}}

\newcommand{\Id}{{\operatorname{Id}}}

\newcommand{\Spec}{\operatorname{Spec}}

\newcommand{\depth}{\operatorname{depth}}
\newcommand{\codepth}{\operatorname{codepth}}
\newcommand{\codim}{\operatorname{codim}}

\newcommand{\hdim}{\operatorname{hdim}}
\newcommand{\supp}{\operatorname{supp}}
\newcommand{\bbP}{{\mathbb P}}
\newcommand{\cA}{{\mathcal A}}
\newcommand{\cB}{{\mathcal B}}

\newcommand{\cE}{{\mathcal E}}
\newcommand{\cF}{{\mathcal F}}

\newcommand{\calH}{{\mathcal H}}
\newcommand{\cI}{{\mathcal I}}
\newcommand{\cJ}{{\mathcal J}}
\newcommand{\calL}{{\mathcal L}}
\newcommand{\cK}{{\mathcal K}}
\newcommand{\cM}{{\mathcal M}}
\newcommand{\cN}{{\mathcal N}}
\newcommand{\cO}{{\mathcal O}}
\newcommand{\cP}{{\mathcal P}}

\newcommand{\bR}{{\mathbf R}}
\newcommand{\bL}{{\mathbf L}}

\newcommand{\cplx}[1]{{{\mathcal #1}^{\scriptscriptstyle\bullet}}}

\newcommand{\dcplx}[1]{{{\mathcal #1}^{\scriptscriptstyle\bullet\vee}}}
\newcommand{\dSHom}[1]{{\mathcal{H}om_{#1}^{\scriptscriptstyle\bullet}}}

\newcommand{\marginnote}[1]{\ifthenelse{\isodd{\thepage}}{\normalmarginpar}
{\reversemarginpar}\marginpar{\fbox{\parbox{24mm}{\sloppy\footnotesize #1}}}}

\newcommand{\rest}[2]{{#1}_{\vert #2}}
\newcommand{\iso}{{\,\stackrel {\textstyle\sim}{\to}\,}}
\newcommand{\cono}{\operatorname{Cone}}

 \newtheorem{thm}{Theorem}[section]
 \newtheorem*{thm*}{Theorem}
 \newtheorem{cor}[thm]{Corollary}
 \newtheorem{lem}[thm]{Lemma}
 \newtheorem{prop}[thm]{Proposition}
 
 \theoremstyle{definition}
 \newtheorem{defin}[thm]{Definition}
 \newenvironment{defn}{\begin{defin}}{\hfill\hspace{1pt}$\triangle$\end{defin}}
 \theoremstyle{remark}
 \newtheorem{rema}[thm]{Remark}
 \newtheorem{exe}[thm]{Example}
 \newenvironment{ex}{\begin{exe}}{\hfill\hspace{1pt}$\triangle$\end{exe}}
\newenvironment{rem}{\begin{rema}}{\hfill\hspace{1pt}$\triangle$\end{rema}}

\numberwithin{equation}{section}
\begin{document}
\title{FOURIER MUKAI TRANSFORMS FOR GORENSTEIN SCHEMES}
\author[D. Hern\'andez Ruip\'erez]{Daniel Hern\'andez Ruip\'erez}
\email{ruiperez@usal.es}
\address{Departamento de Matem\'aticas and Instituto Universitario de F\'{\i}sica Fundamental y Matem\'aticas
(IUFFYM), Universidad de Salamanca, Plaza
de la Merced 1-4, 37008 Salamanca, Spain}
\author[A.C. L\'opez Mart\'{\i}n]{Ana Cristina L\'opez Mart\'{\i}n}
\email{anacris@usal.es}
\author[F. Sancho de Salas]{Fernando Sancho de Salas}
\email{fsancho@usal.es}
\address{Departamento de Matem\'aticas, Universidad de Salamanca, Plaza
de la Merced 1-4, 37008 Salamanca, Spain}
\date{\today}
\thanks {Work supported by research projects BFM2003-00097 (DGI)
%``Transformadas Geom\'etricas Integrales y Aplicaciones''
and SA114/04 (JCYL)}
\begin{abstract}
We extend to singular schemes with Gorenstein
singularities or fibered in schemes of that kind Bondal and Orlov's criterion for an integral functor to be fully faithful. We also prove that the original condition of characteristic zero cannot be removed by providing a counterexample in positive characteristic. We contemplate a criterion for equivalence as well.  In addition, we prove that for locally projective Gorenstein morphisms, a relative
 integral functor is fully faithful if and only if its restriction to each fibre is also  fully faithful. These
 results imply the invertibility of the usual relative Fourier-Mukai transform for an elliptic fibration as a
 direct corollary.
\end{abstract}
\maketitle
%%%%%%%%%%

{\small \tableofcontents }

\section*{Introduction}

Since its introduction by Mukai  \cite{Muk81}, the theory of
integral functors and Fourier-Mukai transforms have been important tools
in the study of the geometry of varieties and moduli spaces. At
the first moment, integral functors were used mainly in connection
with moduli spaces of sheaves and bundles, and provided new
insights in the theory of Picard bundles on abelian varieties and
in the theory of stable sheaves on abelian or K3 surfaces
\cite{Muk87a,BBH97a,BBH98a,Bri99}. In the relative version
\cite{Muk87b,BBHM98,Bri98,BBHM02,HMP02,BrM02,Ca2,BuKr04a} they
have been also used in mirror symmetry and to produce new
instances of stable sheaves on elliptic surfaces or elliptic
Calabi-Yau threefolds. The reason is that the theory of integral
functors is  behind  the spectral data constructions
\cite{FMW99,ACHY01,Ho01}; the irruption of the derived
categories in string theory  caused by homological mirror symmetry
brought then a new interest to derived categories and integral
functors (see \cite{AnHR05,As04a} for recent surveys of the
subject and references therein).

Aside from their interest in Physics, derived categories are
important geometric invariants of algebraic varieties. Much work
is being done in this direction, particularly in the
characterisation of all the algebraic varieties sharing the same
derived category (also known as Mukai partners).

There are classic results like the theorem of Bondal and Orlov
\cite{BO01}  which says that if $X$ is a smooth projective variety
whose canonical divisor is either  ample or anti-ample,
then $X$ can be reconstructed from its
derived category. Mukai proved \cite{Muk81} that there exist non
isomorphic abelian varieties and non isomorphic K3 surfaces having
equivalent derived categories. Orlov \cite{Or97}  proved that two
complex K3 surfaces have equivalent derived categories if and only
if the transcendental lattices of their cohomology spaces are
Hodge-isometric, a result now called the derived Torelli theorem for K3
surfaces. After Mukai's work the problem of finding Fourier-Mukai
partners has been contemplated by many people. Among them, we can
cite Bridgeland-Maciocia \cite{BrM01} and  Kawamata \cite{Kaw02b};
they have proved  that if $X$ is a smooth projective surface, then
there is a finite number of surfaces $Y$ (up to isomorphism)
whose derived category is equivalent to the derived category of
$X$. Kawamata proved that if $X$ and $Y$ are smooth projective
varieties with equivalent derived categories, then $n=\dim X=\dim
Y$ and if moreover $\kappa(X)=n$ (that is, $X$ is of general
type), then there exist birational morphisms $f\colon Z\to X$,
$g\colon Z\to Y$ such that $f^\ast K_X\sim g^\ast K_Y$
(i.e.~$D$-equivalence implies $K$-equivalence) \cite{Kaw02b}.
Other important contributions are owed to Bridgeland \cite{Bri02},
who proved that two crepant resolutions of a projective threefold
with terminal singularities have equivalent derived categories;
therefore,  two birational Calabi-Yau threefolds have equivalent
derived categories. The proof is based on a careful study of the
behaviour of flips and flops under certain integral functors and the construction of the moduli space of perverse point sheaves.

All these results support the belief that derived categories and
integral functors could be most useful in the understanding of the
minimal model problem in higher dimensions. And this suggests that
the knowledge of both the derived categories and the properties of
integral functors for singular varieties could be of great
relevance.

However, very little attention has been paid so far to singular
varieties in the Fourier-Mukai literature. One of the reasons may
be the fact that the fundamental results on integral functors
are not easily generalised to the singular situation,  because
they rely deeply on properties inherent to smoothness.

We would like to mention two of the most important. One is Orlov's
representation theorem \cite{Or97} according to   if $X$ and $Y$
are smooth projective varieties, any (exact) fully faithful
functor between their derived categories is an integral functor.
Particularly, any (exact) equivalence between their derived
categories is an integral functor (integral functors that are
equivalences are also known as Fourier-Mukai functors). Another is
Bondal and Orlov's characterisation of those integral functors
between the derived categories of two smooth varieties that are
fully faithful  \cite{BO95}.

Orlov's representation theorem has been generalised by Kawamata
\cite{Kaw02} to the smooth stack associated to a normal projective
variety with only quotient singularities. Therefore
$D$-equivalence also implies $K$-equivalence for those varieties when $\kappa(X)$ is maximal.
In  \cite{VdB04} Van den Bergh proves using non-commutative rings
that Bridgeland's result about flopping contractions can be
extended to quasi-projective varieties with only Gorenstein
terminal singularities. The same result was proved by Chen
\cite{Chen02};  the underlying idea is to embed such a threefold
into a smooth fourfold and then use the essential smoothness. The
author himself notices that his smoothing approach will not work
for most general threefold flops because quotient singularities in
dimension greater or equal to 3 are very rigid. In his paper, some
general properties of the Fourier-Mukai transform on singular
varieties can be found as well as the computation of a spanning
class of the derived category of a normal projective variety with
only isolated singularities. Finally, Kawamata \cite{Kaw02a} has
obtained analogous results for some $\mathbb{Q}$-Gorenstein
threefolds using algebraic stacks.

This paper is divided in two parts. In the first part, we give an
extension of Bondal and Orlov's characterisation of fully faithful
integral functors to proper varieties with (arbitrary) Gorenstein
singularities. This is the precise statement.
\begin{thm*}[Theorem \ref{1:ffcritGo}]   Let $X$ and $Y$ be projective Gorenstein schemes over an algebraically closed field of characteristic zero, and let $\cplx{K}$ be an object in $\cdbc{X\times Y}$  of finite
projective dimension over $X$ and over $Y$. Assume also that $X$
is
%irreducible
integral. Then the functor $\fmf{\cplx{K}}{X}{Y}\colon
\cdbc{X}\to \cdbc{Y}$ is fully faithful if and only if the kernel
$\cplx{K}$ is strongly simple over $X$.
\end{thm*}

One should notice that this Theorem may fail to be true in positive characteristic even in the smooth case. A counterexample is given in Remark \ref{r:charp}.

In the Gorenstein case, strong simplicity (Definition
\ref{1:strgsplcplxGo}) is defined in terms of locally complete
intersection zero cycles instead of the structure sheaves of the
closed points, as it happens in the smooth case. In the latter
situation, our result improves the characterization of fully
faithfulness of Bondal and Orlov.

As in the smooth case, when $X$ is a Gorenstein variety the
skyscraper sheaves $\cO_x$ form a spanning class for the derived
category $\cdbc{X}$.  Nevertheless, due to the fact that one may has an infinite number of non-zero $\Ext^i_X(\cO_x,\cO_x)$ when $x$ is a singular point, this spanning class does not allow to give an effective criterion characterising the fully faithfulness of
integral functors. However, Bridgeland's criterion that characterises when a fully faithful integral functor is an equivalence is also valid in the Gorenstein case. Moreover, since for a Gorenstein variety one has a more natural
spanning class given by the structure sheaves of locally
complete intersection cycles supported on closed points, one
also proves the following alternative result.

\begin{thm*}[Theorem \ref{t:equivalence}] Let $X$, $Y$ and $\cplx{K}$  be as
in the previous theorem with $Y$ connected. A fully faithful integral functor
$\fmf{\cplx{K}}{X}{Y}\colon \cdbc{X}\to \cdbc{Y}$ is an equivalence of categories
 if and only if for every closed point $x\in X$ there exists a  locally complete
 intersection cycle $Z_x$ supported on $x$ such that $\fmf{\cplx K}{X}{Y}(\cO_{Z_x})
 \simeq \fmf{\cplx K}{X}{Y}(\cO_{Z_x})\otimes \omega_Y$.
\end{thm*}

We also derive in the Gorenstein case some geometric consequences
of the existence of Fourier-Mukai functors (Proposition
\ref{p:conseq}) which are analogous to certain well-known
properties of smooth schemes.

The second part of the paper is devoted to relative integral
functors. As already mentioned, relative
Fourier-Mukai transforms have been considered mainly in connection
with elliptic fibrations. And besides some standard functorial
properties, like compatibility with (some) base changes, more
specific results or instances of Fourier-Mukai functors
(equivalences of the derived categories) are known almost only for
abelian schemes \cite{Muk87b} or elliptic fibrations.

We prove a new result that characterises when  a relative integral
functor is fully faithful or an equivalence, and generalises
\cite[Prop.~6.2]{Chen02}:

\begin{thm*}[Theorem \ref{t:relative}]
Let $p\colon X\to S$ and $q\colon Y\to S$ be locally projective
Gorenstein morphisms (the base field is algebraically closed of
characteristic zero). Let $\cplx{K}\in \dbc{X\times_S Y}$ be a
kernel  of finite projective dimension over both $X$ and $Y$. The
relative integral functor $\fmf{\cplx{K}}{X}{Y}\colon \cdbc{X}\to
\cdbc{Y}$ is fully faithful (respectively an equivalence) if and
only if $\fmf{\bL j_s^\ast\cplx{K}}{X_s}{Y_s}\colon \cdbc{X_s}\to
\cdbc{Y_s}$ is fully faithful (respectively an equivalence) for
every closed point $s\in S$, where $j_s$ is the immersion of
$X_s\times Y_s$ into $X\times_S Y$.
\end{thm*}

Though this result is probably true in greater generality, our
proof needs the Gorenstein condition in an essential way. The
above theorem, together with the characterisation of fully
faithful integral functors and of Fourier-Mukai functors in the
absolute Gorenstein case (Theorems \ref{1:ffcritGo} and
\ref{t:equivalence}) gives a criterion to ascertain when a
relative integral functor between the derived categories of the
total spaces of two Gorenstein fibrations is an equivalence. We
expect that this theorem could be applied to very general
situations. As a first application we give here a very simple and
short proof of  the invertibility result for elliptic fibrations:
\begin{thm*}[Proposition \ref{p:poincequiv}]
Let $S$ be an algebraic scheme over an algebraically closed field of characteristic zero, $X \to S$ an elliptic fibration
with integral fibres and a section, $\hat X\to S$ the dual
fibration and $\cP$ the relative Poincar\'e sheaf on
$X\times_S\hat X$. The relative integral functor
$$\fmf{\cP}{X}{\hat X}\colon \cdbc{X}\to \cdbc{\hat X}$$  is an equivalence of categories.
\end{thm*}

This result has been proved elsewhere in different ways.  When the
total spaces of the fibrations involved are smooth the theorem can
be proved, even if the fibres are singular, by considering the
relative integral functor as an absolute one (defined by the
direct image of the relative Poincar\'e to the direct product) and
then applying the known criteria in the smooth case
\cite{Bri98,BrM02,BBHJ} (see also \cite{BBHM98}). When the total
spaces are singular, there is a proof in \cite{BuKr04,BuKr04a}
that follows a completely different path and is much longer than
ours.

\medskip

In this paper, scheme means algebraic scheme (that is, a scheme of finite type) over an algebraically closed field $k$. By a Gorenstein morphism, we understand a flat morphism of schemes whose fibres are Gorenstein. For any scheme $X$ we denote by $\catD(X)$ the
derived category of complexes of $\cO_X$-modules with
quasi-coherent cohomology sheaves. This is the essential image of
the derived category of quasi-coherent sheaves in the derived
category of all $\cO_X$-modules. Analogously $\dpc{X}$, $\dmc{X}$
and $\dbc{X}$ will denote the derived categories of complexes
which are respectively bounded below, bounded above and bounded on
both sides, and have quasi-coherent cohomology sheaves. The
subscript $c$ will refer to the corresponding subcategories of
complexes with coherent cohomology sheaves.

\subsection*{Acknowledgements} We would like to thank Claudio Bartocci, Ugo Bruzzo
and Tom Bridgeland for discussions, indications and ideas.  We also thank Alessio Corti
 and Constantin Teleman for organising a discussion session about Fourier-Mukai for
 singular varieties during the ``Coherent Sheaves and Mirror Symmetry'' workshop held
 in Cambridge (May 2005) and to the participants in the session, specially Igor Burban,
 Alexander Kuznetsov, Dimitri Orlov, Tony Pantev and Miles Reid, for useful comments and suggestions.
Finally we owe our gratitude to the authors of the forthcoming book \cite{BBHJ} for sharing with us their notes and to the anonymous referee for comments and suggestions which helped us to improve the manuscript.

\section{Fourier-Mukai transform on Gorenstein schemes}

\subsection{Preliminary results}

We first recall some basic formulas which will be used in the rest of the paper.

If $X$ is a scheme, there is  a functorial isomorphism (in the derived
category)
\begin{equation}
\label{tens1}
 \bR \dSHom{\cO_X} (\cplx{F},\bR \dSHom{\cO_X} (\cplx{E},\cplx{H})) \iso
 \bR \dSHom{\cO_X} (\cplx{F}\lotimes \cplx{E}, \cplx{H})
 \end{equation}
  where
$\cplx{F}$,  $\cplx{E}$ are in $\dmc{X}$, $\cplx{H}$ is in $\dpc{X}$, and all have coherent cohomology (\cite{Hart66}).  One also has a functorial isomorphism
\begin{equation}
\label{tens2}
 \bR \dSHom{\cO_X} (\cplx{F},\cplx{E}) \lotimes \cplx{H} \iso \bR \dSHom{\cO_X}
 (\cplx{F},\cplx{E}\lotimes \cplx{H})
 \end{equation}
where $\cplx{F}$ is a bounded complex of $\cO_X$-modules with coherent
cohomology and either $\cplx F$ or $\cplx H$ is of finite
homological dimension  (i.e.~ locally isomorphic to a bounded complex of locally free sheaves of finite rank). The usual proof (see \cite{Hart66} or
\cite{BBHJ}) requires that $\cplx{H}$ is of finite homological
dimension; however, it still works when both members are defined. If we denote by $\cplx{F}^\vee = \bR \dSHom{\cO_X} (\cplx{F},\cO_X)$ the dual in the derived category, \eqref{tens2} implies that
\begin{equation}
\label{tens3}
\cplx{F}^\vee \lotimes \cplx{H} \iso \bR \dSHom{\cO_X}
 (\cplx{F}, \cplx{H})\,.
 \end{equation}
Nevertheless this formula may fail to be true when neither $\cplx
F$ nor $\cplx H$ have finite homological dimension as the
following example shows.

\begin{ex} Let $X$ be a Gorenstein scheme of dimension $n$ over a
field $k$. Let $x\in X$ be a singular point and let $\cF$ be
any $\cO_X$-module. Since $\cO_x^{\scriptscriptstyle\vee}\simeq
\cO_x[-n]$, if one had
$$\cO_x^{\scriptscriptstyle\vee}\lotimes\cF\simeq
\bR \dSHom{\cO_X} (\cO_x,\cF)\, ,$$ then one would have $
\operatorname{Tor}_{n-i}(\cO_x,\cF)\simeq\Ext^i(\cO_x,\cF)$ for
every $i\in \mathbb{Z}$. It follows that $\Ext^i(\cO_x,\cF)=0$ for
all $i>n$ and every $\cO_X$-module $\cF$ and this is impossible
because $\cO_x$ is not of finite homological dimension.
\end{ex}

The formula \eqref{tens3} implies that if $f\colon X \to Y$ is a
morphism, there is an isomorphism
\begin{equation}
\label{basechangedual} \bL f^\ast (\dcplx{F}) \simeq (\bL f^\ast \cplx F)^{\scriptscriptstyle\vee}
\end{equation}
if either $\cplx F$ is of finite homological dimension or $f$ is
of finite Tor-dimension (in this paper we shall only need to consider the case
when $f$ is flat or is a regular closed immersion).
% $\cO_X$ is of
%finite homological dimension as a module over $\cO_Y$ (i.e., it is
%locally isomorphic to a bounded complex of sheaves which are flat
%over $\cO_Y$).

Some other formulas will be useful. When $X$ is a  Gorenstein
scheme, every object  $\cplx F$ in $\cdbc{X}$ is reflexive, that is, one has an isomorphism in the derived category \cite[1.17]{Or03}:
\begin{equation}
\label{ddual} \cplx{F}\simeq(\dcplx{F})^{\scriptscriptstyle\vee}\,.
\end{equation}
Then, one has
\begin{equation}\label{homdual}
\Hom_{\catD(X)}(\cplx{H},\cplx{F})\simeq
\Hom_{\catD(X)}(\dcplx{F},\dcplx{H})
\end{equation}
for every bounded complex $\cplx{F}$ in $\cdbc{X}$ and any complex
$\cplx{H}$.

Moreover, if $X$ is a zero dimensional Gorenstein scheme, the sheaf $\cO_X$ is injective so that
\begin{equation} \label{e:zeroddual}
\cplx F^\ast\simeq \cplx F^\vee \quad \text{and} \quad \calH^i(\cplx F^\vee)\simeq (\calH^{-i}(\cplx F))^\ast
\end{equation}
 for every object  $\cplx F$ in $\cdbc{X}$, where $\cplx F^\ast = \dSHom{\cO_X} (\cplx F,\cO_X)$ is the ordinary dual.

Let $f\colon X\to Y$ be a proper morphism of schemes. The relative
Grothendieck duality states the existence of a functorial
isomorphism in the derived category
\begin{equation}
\label{1:duality} \bR\dSHom{\cO_Y}(\bR f_{\ast}\cplx{F},\cplx{G})
\simeq \bR f_\ast \bR\dSHom{\cO_{X}} (\cplx{F}, f^!\cplx{G})\,.
\end{equation} for $\cplx G$ in $\catD(Y)$ and $\cplx F$ in
$\catD(X)$ (see for instance \cite{Nee96}). By applying the
derived functor of the global section functor, we obtain the
\emph{global duality formula}
\begin{equation} \label{1:adj1}
\Hom_{\catD(Y)}(\bR f_{\ast}\cplx{F},\cplx{G}) \simeq
\Hom_{\catD(X)} (\cplx{F}, f^!\cplx{G})\,.
\end{equation}
In other words, the direct image $\bR f_{\ast}\colon \catD(X)\to
\catD(Y)$ has a right adjoint $ f^!\colon \catD(Y) \to \catD(X)$.

There is a natural map $f^* \cplx{G}\lotimes f^!\cO_Y\to
f^!\cplx{G}$, which is an isomorphism when either $\cplx{G}$ has finite homological dimension or $\cplx{G}$ is reflexive and  $f^!\cO_Y$ has finite homological dimension.

When $f$ is a Gorenstein morphism of relative dimension $n$, the object $f^!\cO_Y$  reduces to an invertible sheaf $\omega_f$, called \emph{the relative dualizing sheaf}, located at the place $-n$, $f^!\cO_Y\simeq \omega_f[n]$.

Grothendieck duality is compatible with base-change. We state this result for simplicity only when $f$ is Gorenstein. In this case,  since $f$ is flat, base-change compatibility means that if $g\colon Z\to Y$
is a morphism and $f_Z\colon Z\times_Y X\to Z$ is the induced
morphism, then the relative dualizing sheaf for $f_Z$ is
$\omega_{f_Z}=g_X^\ast \omega_f$ where $g_X\colon Z\times_Y X\to
X$ is the projection.

As it is customary, when $f$ is the projection onto a point, we
denote the dualizing sheaf by $\omega_X$.

\subsection{Complexes of relative finite projective dimension}

In this subsection we shall prove a weaker version of \eqref{tens2} in some cases.
%Besides the usual relative notion of finite homological dimension, we need to consider the notion  of finite projective dimension.
\begin{lem} \label{locallyperfect} Let $\cplx{E}$ be an object in $\cdbc{X}$. The
following conditions are equivalent:
\begin{enumerate}\item $\cplx{E}$ is of finite homological
dimension. \item $\cplx E\lotimes \cplx G$  is an object of
$\dbc{X}$ for every $\cplx G$ in $\dbc{X}$. \item
$\bR\dSHom{\cO_{X}} (\cplx{E}, \cplx{G})$ is in $\dbc{X}$ for
every $\cplx{G}$ in $\dbc{X}$.
\end{enumerate}
\end{lem}

\begin{proof} Since $X$ is noetherian, the three conditions are local
so that we can assume that $X$ is affine. It is clear that (1)
implies (2) and (3). Now let us see that (3) implies (1). Let us
consider a quasi-isomorphism $\cplx L \to \cplx E$ where $\cplx L$
is a bounded above complex of finite free modules. If $\cK^n$ is
the kernel of the differential $\calL^n \to \calL^{n+1}$,  then
for $n$ small enough the truncated complex $\cK^n \to \calL^n \to
\dots $ is still quasi-isomorphic to $\cplx E$ because $\cplx E$
is an object of $\cdbc{X}$. Let $x$ be a point and $\cO_x$ its
residual field. Since $\bR\dSHom{\cO_X} (\cplx E, \cO_x)$ has
bounded homology, one also has that
$\Ext^1_{\cO_X}(\cK^n,\cO_x)=0$ for $n$ small enough. For such $n$
the module $\cK^n$ is free in a neighbourhood of $x$ and one
concludes. To prove that (2) implies (1), one proceeds analogously
replacing $\Ext^1$ by $\Tor_1$.
\end{proof}
This lemma suggests the following definition.

\begin{defn} Let $f\colon X\to Y$ be a morphism of schemes.
An object $\cplx E$ in $\catD(X)$ is said to be of \emph{finite
homological dimension over $Y$} (resp.~of \emph{finite projective
dimension over $Y$}),  if  $\cplx E\lotimes \bL f^\ast \cplx G$
(resp.~ $\bR\dSHom{\cO_{X}} (\cplx{E}, f^!\cplx{G})$),  is in
$\dbc{X}$ for any $\cplx{G}$ in $\dbc{Y}$.
\end{defn}

These notions are similar (though weaker) to the notions of finite Tor-amplitude and finite Ext-amplitude considered in \cite{Kuz05}.
%If $\cO_X$ is of finite homological dimension as a sheaf of $\cO_Y$-modules,
%then $f$ is of finite homological dimension in the sense of this definition.

In the absolute case (i.e.~when $f$ is the identity), finite
projective dimension is equivalent to finite homological dimension
by the previous lemma. To characterise complexes of finite projective
dimension over $Y$ when $f$ is projective, we shall need the
following result (c.f.~\cite[Lem.~2.13]{Or97}).
\begin{lem} \label{l:beilinson}
 Let $A$ be a noetherian ring, $f\colon X \to Y=\Spec A$ a projective morphism and $\cO_X(1)$ a relatively very
 ample line bundle.
 \begin{enumerate}
 \item Let $\cplx M$ be an object of  $\dmc{X}$. Then $\cplx M=0$ (resp. is an object of $\dbc{X}$) if and only if $\bR f_\ast (\cplx M(r))=0$
 (resp. is an object of $\dbc{Y}$) for
 every integer $r$.
\item Let $g\colon \cplx M\to\cplx N$ be a morphism in $\dmc{X}$. Then $g$ is an isomorphism if and only if the
induced morphism $\bR f_\ast (\cplx M(r))\to \bR f_\ast (\cplx N(r))$ is an isomorphism in $\dmc{Y}$ for every
integer $r$.
\end{enumerate}
(As it is usual, we set $\cplx M(r)=\cplx M\otimes \cO_X(r)$.)
\end{lem}
\begin{proof} Let $i\colon X\hookrightarrow \bbP_A^N$ be the closed immersion of $A$-schemes defined by  $\cO_X(1)$.
Since $\cplx M=0$ if and only if $i_\ast\cplx M=0$ and $\cplx M$ has bounded cohomology  if and only if $i_\ast\cplx M$ has bounded cohomology as well, we can assume that $X=\bbP_A^N$. Now one has an exact
sequence (Beilinson's resolution of the diagonal)
$$
0\to \cE_N\to \cdots \to
\cE_1\to \cE_0 \to\cO_{\Delta}\to 0
$$
where $\cE_j=\pi_1^*\cO_{\bbP_A^N/A}(-j)\otimes \pi_2^\ast \Omega^j_{\bbP_A^N/A}(j)$,  $\pi_1$ and $\pi_2$ being  the
projections of  $\bbP_A^N\times \bbP_A^N$ onto its factors. Then $\cO_\Delta$ is an object of the smallest
triangulated subcategory of $\dbc{\bbP_A^N\times \bbP_A^N}$ that contains the sheaves $\cE_j$ for $0\le j \le N$.
Since $F(\cplx F)=\bR \pi_{2\ast}(\pi_1^\ast(\cplx M)\lotimes \cplx F)$ is an exact functor
$\dbc{\bbP_A^N\times \bbP_A^N} \to \dmc{\bbP_A^N}$, $\cplx M \simeq F(\cO_\Delta)$ is an object of the smallest
triangulated category generated by the objects $F(\cE_j)$ for $0\le j\le N$. Thus to prove (1)
we have only to see that $F(\cE_j)=0$ (resp. have bounded homology) for all $0\le j\le N$. This follows because we have
$$
F(\cE_j)\simeq \bR \pi_{2\ast}(\pi_1^\ast(\cplx M(-j)))\lotimes \Omega^j_{\bbP_A^N/A}(j)
\simeq f^\ast \bR f_\ast (\cplx M(-j)))\otimes \Omega^j_{\bbP_A^N/A}(j)
$$
by the projection formula \cite[Prop.~5.3]{Nee96} and flat base-change.

By applying the first statement to the cone of $g$, the second statement follows.
\end{proof}

One can also easily prove that $\cplx M=0$ if and only if $\bR f_\ast (\cplx M(r))=0$ for all $r$ by using the spectral sequence $R^p f_\ast (\calH^q(\cplx M(r)))\implies R^{p+q} f_\ast (\cplx M(r))$.

\begin{lem} \label{fpdchar}  Let $f\colon X\to Y$ be a proper morphism and $\cplx{E}$ an object of $\cdbc{X}$. If $\cplx E$ is either of finite projective dimension or of finite homological dimension over $Y$, then
$\bR f_\ast \cplx E$ is of finite homological dimension.
\end{lem}
\begin{proof} The duality isomorphism \eqref{1:duality}  together with Lemma \ref{locallyperfect} imply that
$\bR f_\ast \cplx E$ is of finite homological dimension when  $\cplx E$ is of finite
projective dimension over $Y$. If  $\cplx E$ is of finite
homological dimension over $Y$, we use the same lemma and the
projection formula.
\end{proof}

\begin{prop}\label{fpdchar2}  Let $f\colon X\to Y$ be a projective morphism and $\cplx{E}$ an object of $\cdbc{X}$.
 The following conditions are equivalent:
\begin{enumerate}
\item   $\cplx E$ is of finite projective dimension over $Y$.
\item $\bR f_\ast
(\cplx E(r))$ is of finite homological dimension for every
integer $r$.
\item $\cplx E$ is of finite homological dimension over $Y$.
\end{enumerate}
Thus, if $f$ is locally projective, $\cplx E$ is of finite projective dimension over $Y$ if and only if it is of  finite homological dimension over $Y$.
\end{prop}
\begin{proof}  If $\cplx E$ is of finite projective dimension (resp.~of finite homological dimension) over $Y$, so is $\cplx E(r)$
 for every $r$, and then  $\bR f_\ast
(\cplx E(r))$ is of finite homological dimension by Lemma \ref{fpdchar}. Assume that (2)
is satisfied.  Then (1) is a consequence of the duality
isomorphism $\bR f_\ast (\bR\dSHom{\cO_{X}} (\cplx{E},
f^!\cplx{G})(r)) \simeq \bR\dSHom{\cO_Y}(\bR
f_{\ast}(\cplx{E}(-r)),\cplx{G})$ and Lemma \ref{locallyperfect},
whilst (3) follows from the same lemma and the projection formula.
\end{proof}

\begin{cor} \label{fpddual} Let $f\colon X\to Y$ be a projective morphism and $\cplx{E}$ an object of $\cdbc{X}$. If $\cplx E$ is of finite projective dimension over $Y$, then
$\bR \dSHom{\cO_X} (\cplx E, f^!\cO_Y)$  is also of finite
projective dimension over $Y$. In particular, if $f$ is
Gorenstein, $\cplx{E}^\vee$ is of finite projective dimension over
$Y$.
\end{cor}
\begin{proof}  Let us write $\cplx N=\bR\dSHom{\cO_{X}} (\cplx{E}, f^!\cO_Y)$. By Proposition \ref{fpdchar2}, it suffices to see that
$\bR f_\ast(\cplx N(r))$ is of finite homological dimension for every $r$. This follows again by Proposition \ref{fpdchar2}, due to the isomorphism
$\bR f_\ast(\cplx N(r))\simeq  [\bR f_\ast (\cplx E(-r))]^\vee$.
\end{proof}

\begin{prop}
\label{dualGore} Let $f\colon X\to Y$ be a locally projective
Gorenstein morphism of  schemes and $\cplx{E}$ an object of
$\catD_c(X)$ of finite projective dimension over $Y$. One has
$$
\dcplx{E}\lotimes  f^\ast\cplx{G}\otimes \omega_f [n] \simeq \bR
\dSHom{\cO_{X}} (\cplx{E}, f^!\cplx{G})
$$
for $\cplx{G}$ in $\cdbc{Y}$.
%In particular, $\dcplx{E}$ is of
%finite homological dimension over $Y$.
Moreover, if $Y$ is
Gorenstein, then
$$
\dcplx{E}\lotimes  f^\ast\cplx{G} \simeq \bR \dSHom{\cO_{X}}
(\cplx{E}, f^\ast\cplx{G})\,.
$$
\end{prop}
\begin{proof} One has natural morphisms
\begin{equation} \label{eqiso}
\begin{aligned} \bR\SHom_{\cO_X}(\cplx{E},\cO_X)\lotimes
f^\ast\cplx{G}\otimes \omega_f[n]&\to \bR\SHom_{\cO_X}(\cplx{E},
f^\ast\cplx{G}\otimes \omega_f[n]) \\ &\to
\bR\SHom_{\cO_X}(\cplx{E}, f^! \cplx{G})\,.
\end{aligned}
\end{equation}
 We have to prove that the composition is an isomorphism. This is
a local question on $Y$, so that we can assume that  $Y=\Spec A$.

By Lemma \ref{l:beilinson} we have to prove that
the induced morphism
\begin{equation}
\bR f_{\ast}(\bR\SHom_{\cO_X}(\cplx{E},\cO_X)\lotimes
f^\ast\cplx{G}\otimes\omega_f[n]\otimes \cO_{X}(r)) \to \bR
f_{\ast}(\bR\SHom_{\cO_X}(\cplx{E}, f^! \cplx{G})\otimes
\cO_{X}(r))
\end{equation}
is an isomorphism in $\dmc{Y}$ for any integer $r$.
 The first member is isomorphic to
\[
\bR\SHom_{\cO_Y}(\bR f_\ast(\cplx{E}(-r) ),\cO_Y)\lotimes \cplx{G}
\]
by the projection formula and relative duality; the second one is
isomorphic to
\begin{equation} \label{e:bounded}
\bR f_\ast (\bR\SHom_{\cO_X}(\cplx{E}, f^! \cplx{G})\otimes
\cO_X(r)) \simeq \bR\SHom_{\cO_Y}(\bR f_\ast(\cplx{E}(-r)),
\cplx{G})\,.
\end{equation}
Thus, we have to prove that the natural morphism
\begin{equation}
\bR\SHom_{\cO_Y}(\bR f_\ast \cplx{E}(-r),\cO_Y)\lotimes \cplx{G} \to
\bR\SHom_{\cO_Y}(\bR f_\ast \cplx{E}(-r), \cplx{G})\,,  \label{eqiso2}
\end{equation}
is an isomorphism. Since $\bR f_\ast \cplx{E}(-r)$ is of finite homological dimension by Proposition \ref{fpdchar2}, one concludes by \eqref{tens2}.
\end{proof}

\subsection{Depth and local properties of Cohen-Macaulay and Gorenstein schemes}

Here we state some preliminary results about depth on singular
schemes and local properties of Cohen-Macaulay and Gorenstein
schemes. We  first recall a local property of Cohen-Macaulay
schemes.

\begin{lem}\label{l:zerociclo} \cite[Prop.~6.2.4]{Rob98} Let $A$ be a noetherian local ring.
$A$ is Cohen-Macaulay if and only if there is an ideal $I$ of $A$
with $\dim A/I=0$ and such that $A/I$ has finite homological
dimension.
\end{lem}
\begin{proof} Let $n$ be the dimension of the ring $A$.
If $A$ is Cohen-Macaulay, $\depth(A)=n$. Then there is a regular
sequence $(a_1,\hdots,a_n)$ in $A$ and taking $I=(a_1, \hdots,
a_n)$ we conclude. Conversely, if $I$ is an ideal satisfying $\dim
A/I=0$ and $\hdim(A/I)=s<\infty$, then the Auslander and Buchsbaum's
formula $\depth(A/I)+ \hdim(A/I)=\depth(A)$
\cite[Thm.~19.1]{Mat86} proves that $s=\depth(A)$; then $s\le n$.
Moreover, if $\cplx{M}$ is a free resolution of $A/I$ of length
$s$, one has $s\geq n$ by the intersection theorem
\cite[6.2.2.]{Rob98}. Thus $A$ is Cohen-Macaulay.
\end{proof}

Let $\cF$ be a coherent sheaf on a scheme $X$
of dimension $n$. We write $n_x$ for the dimension of the local ring $\cO_{X,x}$ of $X$
at a point $x\in X$ and $\cF_x$ for the stalk of $\cF$ at $x$.  $\cF_x$ is a $\cO_{X,x}$-module.
The integer number $\codepth(\cF_x)=n_x-\depth(\cF_x)$ is called the codepth of $\cF$ at $x$.
For any integer $m\in\mathbb{Z}$, the \emph{$m$-th
singularity set} of $\cF$ is defined to be
$$
S_m(\cF)= \{x\in X
\;|\; \codepth(\cF_x)\geq n-m\}\, .
$$
Then, if $X$ is equidimensional, a closed point $x$ is in
$S_m(\cF)$ if and only if $\depth(\cF_x)\leq m$.  If $x$ is a
point of $X$ (not necessarily closed) the zero cycles $Z_x$ of
$\Spec \cO_{X,x}$ supported on the closed point $x$ of $\Spec
\cO_{X,x}$ will be called \emph{zero cycles (of $X$) supported on
$x$} by abuse of language.

Since $\depth(\cF_x)$ is the first integer $i$ such that either
\begin{itemize}
\item $\Ext^i(\cO_x,\cF)\neq 0$ or
\item $H^i_x(\Spec\cO_{X,x},\cF_x)\neq 0$ or
\item $\Ext^i(\cO_Z,\cF_x)\neq 0$ for some zero cycle $Z$ supported on $x$ or
\item $\Ext^i(\cO_Z,\cF_x)\neq 0$ for every zero cycle $Z$ supported on $x$
\end{itemize}
(see for instance \cite{Hart67}), we have alternative descriptions of $S_m(\cF)$:
\begin{equation} \label{e:descrip}
\begin{aligned}
S_m(\cF) &= \{x\in X
\;|\; H^i_x(\Spec\cO_{X,x},\cF_x)\neq 0 \text{ for some } i\leq m+n_x-n\} \\
& =  \{x\in X
\;|\; \Ext^i(\cO_Z,\cF_x)\neq 0 \text{ for some } i\leq m+n_x-n \text{ and some} \\
& \qquad \text{zero cycle $Z$ supported on $x$}\} \\
& =  \{x\in X
\;|\; \Ext^i(\cO_Z,\cF_x)\neq 0 \text{ for some } i\leq m+n_x-n \text{ and any} \\
& \qquad \text{zero cycle $Z$ supported on $x$}\}
\end{aligned}
\end{equation}

\begin{lem} \label{l:sm} If $X$ is smooth, then the $m$-th singularity set of $\cF$ can be described as
$$
 S_m(\cF)=\cup_{p\ge n-m} \{ x\in X \;|\; L_p j_x^\ast \cF\neq
0\}\,,
$$
where $j_x$ is the immersion of the point $x$.
\end{lem}
\begin{proof}  Let $x\in X$ be a point and $\cplx{L}$ the Koszul complex associated locally to a
regular sequence of generators of the maximal ideal of
$\cO_{X,x}$. Since $\cplx{L}^\vee\simeq \cplx{L}[-n_x]$, one has
an isomorphism $\Ext^{i}(\cO_{x},\cF_x)\simeq L_{n_x-i} j_{x}^\ast
\cF$ which proves the result.
\end{proof}

In the singular case, this characterization of $S_m(\cF)$ is not
true. However, there is a similar interpretation for
Cohen-Macaulay schemes as we shall see now. By Lemma
\ref{l:zerociclo}, if $X$ is Cohen-Macaulay, for every point $x$
there exist zero cycles supported on $x$ defined locally by a
regular sequence; we refer to them as  \emph{locally complete
intersection} or \emph{l.c.i.}~cycles. If $Z\hookrightarrow X$ is
such a l.c.i.~cycle, by the Koszul complex theory the structure
sheaf $\cO_Z$ has \emph{finite homological dimension} as an
$\cO_X$-module.

We denote by $j_Z$ the immersion of $Z$ in $X$. Recall that for
every object $\cplx{K}$ in $\dbc{X}$, $L_i j_{Z}^\ast \cplx{K}$
denotes the cohomology sheaf $\calH^{-i}(j_{Z}^\ast\cplx{L})$
where $\cplx{L}$ is a bounded above complex of locally free
sheaves quasi-isomorphic to $\cplx{K}$.

\begin{lem}  \label{l:SetSing} If $X$ is Cohen-Macaulay, then the $m$-th singularity set $S_m(\cF)$ can be described as
\begin{align*}
S_m(\cF) &= \{x\in X \;|\; \text{there is an integer }  i\geq
n-m\text{ with } L_i j_{Z_x}^\ast \cF\neq 0  \\
& \qquad \text{ for any l.c.i
zero cycle $Z_x$ supported on $x$} \}\,.
\end{align*}
\end{lem}
\begin{proof} Let $Z_x$ be a l.c.i.~zero cycle supported on $x$
and $\cplx{L}$ the Koszul complex associated
locally to a regular sequence of generators of the ideal of $Z_x$.
As in the smooth case, we have that $\cplx{L}^\vee\simeq
\cplx{L}[-n_x]$ and then an isomorphism
$\Ext^{i}(\cO_{Z_x},\cF)\simeq L_{n_x-i} j_{Z_x}^\ast \cF$. The
result follows from \eqref{e:descrip}.
\end{proof}

\begin{lem}\label{Simmersion} If $j\colon X\hookrightarrow W$ is a closed immersion
and $\cF$ is a coherent sheaf on $X$, then $S_m(\cF)=S_m(j_\ast
\cF)$.
\end{lem}

\begin{proof} Since  $H^i_x(\Spec\cO_{X,x},\cF_x)
=H^i_x(\Spec\cO_{W,x},(j_\ast\cF)_x)$, the result follows from
\eqref{e:descrip}.
\end{proof}

\begin{prop}\label{p:Sm} Let $X$ be an equidimensional scheme of dimension $n$
and $\cF$ a coherent sheaf on $X$. \begin{enumerate} \item
$S_m(\cF)$ is a closed subscheme of $X$ and $\codim S_m(\cF)\geq
n-m$. \item If $Z$ is an irreducible component of the support of
$\cF$ and $c$ is the codimension of $Z$ in $X$, then $\codim
S_{n-c}(\cF)=c$ and $Z$ is also an irreducible component of
$S_{n-c}(\cF)$.
 \end{enumerate}
\end{prop}

\begin{proof} All questions are local and then, by Lemma
\ref{Simmersion}, we can assume that $X$ is affine and smooth. By
Lemma \ref{l:sm}, $S_m(\cF)=\cup_{p\ge n-m} X_p(\cF)$, where $
X_p(\cF)=\{ x\in X \;|\; L_p j_x^\ast \cF\neq 0\}$. To prove (1),
we have only to see that $X_p(\cF)$ is closed of codimension
greater or equal than $p$. This can be seen by induction on $p$. If
$p=0$, then $X_0(\cF)$ is the support of $\cF$ and the statement
is clear. For $p=1$, $X_1(\cF)$ is the locus of points where $\cF$
is not locally free, which is closed of codimension greater or
equal than 1, since $\cF$ is always free at the generic point. If
$p>1$, let us consider an exact sequence $0\to \cN \to \calL \to
\cF\to 0$ where $\calL$ is free and finitely generated. Then $L_p
j_x^\ast \cF\simeq L_{p-1} j_x^\ast \cN$ so that
$X_p(\cF)=X_{p-1}(\cN)$ which is closed by induction. Moreover, if
$x\in X_p(\cF)$, then $L_p j_x^\ast \cF\neq 0$, so that $p\leq
\dim \cO_{X,x}$ because $\cO_{X,x}$ is a regular ring. It follows
that $\codim X_p(\cF)=\max_{x\in X_p(\cF)} \{\dim \cO_{X,x}\}\geq
p$.

We finally prove (2). By \cite[Thm.~6.5]{Mat86}, the prime ideal
of $Z$ is also a minimal associated prime to $\cF$. Thus, if $x$
is the generic point of $Z$, the maximal ideal of the local ring
$\cO_{X,x}$ is a prime associated to $\cF_x$, and then
$\Hom(\cO_x,\cF_x)\neq 0$. This proves that $x\in S_{n-c}(\cF)$
and then $Z\subseteq S_{n-c}(\cF)$. The result follows.

\end{proof}

\begin{cor}\label{lemafer} Let $X$ be a Cohen-Macaulay scheme  and let $\cF$ be a
coherent $\cO_X$-module.  Let $h\colon Y\hookrightarrow X$ be an irreducible
component of the support of $\cF$ and $c$ the codimension of $Y$
in $X$. There is a non-empty open subset $U$ of $Y$ such that
for any l.c.i.~zero cycle $Z_x$ supported on $x\in U$ one has
\begin{align*}
L_cj_{Z_x}^\ast \cF&\neq 0
\\
L_{c+i}j_{Z_x}^\ast \cF&= 0\,, \quad \text{ for every $i>0$.}
\end{align*}
\end{cor}

\begin{proof} By Lemma \ref{l:SetSing} the locus of the points that verify the conditions is
$U=Y\cap(S_{n-c}(\cF)-S_{n-c-1}(\cF))$, which is open in $Y$ by
Proposition \ref{p:Sm}. Proving that $U$ is not empty is a local
question, and we can then assume that $Y$ is the support of $\cF$.
Now $Y=S_{n-c}(\cF)$ by (2) of Proposition \ref{p:Sm} and
$U=S_{n-c}(\cF)-S_{n-c-1}(\cF)$ is non-empty because
 the codimension of $S_{n-c-1}(\cF)$ in $X$ is greater or equal than $c+1$ again by
Proposition \ref{p:Sm}.
\end{proof}

 The following proposition characterises objects of the derived
category supported on a closed subscheme.
\begin{prop} \emph{\cite[Prop.~1.5]{BO95}}. Let $j\colon Y\hookrightarrow X$ be a closed immersion
of codimension $d$ of irreducible Cohen-Macaulay schemes
and $\cplx{K}$ an object of $\cdbc X$. Assume that
\begin{enumerate}
\item If $x\in X-Y$ is a closed point, then $\bL j_{Z_x}^\ast \cplx{K}=0$ for some
l.c.i.~zero cycle $Z_x$ supported on $x$.
 \item If $x\in Y$ is a closed point, then $L_i
j_{Z_x}^\ast \cplx{K}=0$ for some l.c.i.~zero cycle $Z_x$ supported on $x$ when
either $i<0$ or $i>d$.
\end{enumerate}
Then there is a sheaf $\cK$ on $X$ whose topological support is
contained in $Y$ and such that $\cplx{K}\simeq\cK$ in $\cdbc X$.
Moreover, this topological support coincides with $Y$ unless
$\cplx{K}=0$.
 \label{1:support}
\end{prop}
\begin{proof} Let us
write $\calH^q=\calH^q(\cplx K)$. For every zero cycle $Z_x$ in
$X$ there is a spectral sequence
$$
E_2^{-p,q}=L_p j_{Z_x}^\ast \calH^q \implies
E_\infty^{-p+q}=L_{p-q}j_{Z_x}^\ast \cplx{K}
$$
Let $q_0$ be the maximum of the $q$'s with $\calH^q\neq 0$. If
$x\in \supp(\calH^{q_0})$, one has $j_{Z_x}^\ast \calH^{q_0}\neq
0$ for every l.c.i.~zero cycle $Z_x$ supported on $x$. A nonzero
element in $j_{Z_x}^\ast \calH^{q_0}$ survives up to infinity in
the spectral sequence. Since there is a
l.c.i.~zero cycle $Z_x$ such that $E_\infty^{q}=L_{-q}j_{Z_x}^\ast
\cplx{K}=0$ for every $q>0$ by hypothesis, one has $q_0\le 0$. A similar
argument shows that the topological support of all the sheaves
$\calH^q$ is contained in $Y$: assume that this is not true
and let us consider the maximum $q_1$ of the $q$'s such that
$j_x^\ast\calH^q\neq 0$ for a certain point $x\in X-Y$; then
$j_{Z_x}^\ast\calH^{q_1}\neq 0$ and a nonzero element in
$j_{Z_x}^\ast \calH^{q_1}$ survives up to infinity in the spectral
sequence, which is impossible since $\bL j_{Z_x}^\ast \cplx K=0$.

Let  $q_2\le q_0$ be the minimum of the $q$'s with $\calH^q\neq
0$. We know that $\calH^{q_2}$ is topologically supported on a
closed subset of $Y$. Take a component $Y'\subseteq Y$ of the
support. If $c\ge d$ is the codimension of $Y'$, then there is a
non-empty
 open subset $U$ of $Y'$ such that $L_cj_{Z_x}^\ast \calH^{q_2}\neq 0$ for any
closed point $x\in U$ and any l.c.i.~zero cycle $Z_x$ supported on
$x$, by Corollary \ref{lemafer}. Elements in $L_cj_{Z_x}^\ast
\calH^{q_2}$ would be killed in the spectral sequence by
$L_{p}j_{Z_x}^\ast \calH^{q_2+1}$ with $p\ge c+2$. By Lemma
\ref{l:SetSing}  the set
 $$ \{x\in X \;|\;
L_i j_{Z_x}^\ast \calH^{q_2+1}\neq 0 \text{ for some } i\geq c+2
\text{ and any l.c.i.~cycle } Z_x \}$$ is equal to
$S_{n-(c+2)}(\calH^{q_2+1})$ and then has codimension greater or
equal than $c+2$ by Proposition \ref{p:Sm}. Thus there is a point
$x\in Y'$ such that any nonzero element in $L_cj_{Z_x}^\ast
\calH^{q_2}$ survives up to the infinity in the spectral sequence.
Therefore, $L_{c-q_2}j_{Z_x}^\ast \cplx{K}\neq 0$ for any
l.c.i.~zero cycle $Z_x$ supported on $x$. Thus $c-q_2\le d$ which
leads to $q_2\ge c-d\ge 0$ and then $q_2=q_0=0$. So $\cplx{K}=
\calH^0$ in $\dbc X$ and the topological support of $\cK=\calH^0$
is contained in $Y$. Actually, if $\cplx{K}\neq 0$, then this
support is the whole of $Y$: if this was not true, since $Y$ is
irreducible, the support would have a component $Y'\subset Y$ of
codimension $c>d$ and one could find, reasoning as above, a
non-empty subset $U$ of $Y'$ such that $L_{c}j_{Z_x}^\ast
\cplx{K}\neq 0$ for all $x\in U$ and all l.c.i.~zero cycle $Z_x$
supported on $x$. This would imply that $c\le d$, which is
impossible.
\end{proof}

Taking into account that $\cO_{Z_x}^\vee
=\cO_{Z_x}[-n]$ where $n=\dim X$, Proposition \ref{1:support} may be
reformulated as follows:
\begin{prop}  \label{1:support2} Let $j\colon Y\hookrightarrow X$ be a closed immersion
of codimension $d$ of irreducible Cohen-Macaulay schemes of dimensions $m$ and $n$ respectively,
and let $\cplx{K}$ be an object of $\cdbc X$. Assume that for any closed point $x\in X$ there is a l.c.i.~zero cycle $Z_x$ supported on $x$ such that
$$
 \Hom^i_{\catD(X)}(\cO_{Z_x},\cplx{K})= 0\,,
$$
unless $x\in Y$ and $m\le i\le n$. Then there is a sheaf $\cK$ on
$X$ whose topological support is contained in $Y$ and such that
$\cplx{K}\simeq\cK$ in $\cdbc X$. Moreover, the topological
support is $Y$ unless $\cplx{K}=0$. \qed\end{prop}

%The last result we need is a direct consequence of the
%``intersection theorem'' in commutative algebra \cite[Prop.6.2.4]{Rob98}, \cite{Rob89}, we mentioned earlier.

%\begin{lem}\emph{\cite[Cor.~5.6]{BrM02}}\label{1:intersect} Let $Z$ be an irreducible algebraic
%variety and fix a closed point $x\in Z$. Assume that there is an object $\cplx{E}(x)$ in $\dbc Z$
%such that for any closed point $z\in Z$ and any integer $i$ one has
%$$
%\Hom_{\catD(Z)}^i(\cplx{E}(x),\cO_z)=0\qquad\text{unless $z=x$ and
%$0\le i\le \dim Z$.}
%$$
%Assume also that $\calH^0(\cplx{E}(x))\simeq\sk{x}$. Then $Z$ is
%smooth at $x$ and $\cplx{E}(x)\simeq\sk{x}$. \qed \end{lem}

\subsection{Integral functors}\label{ss:basechange}
Let $X$ and $Y$ be proper schemes. We denote the projections of the direct product
$X\times Y$ to $X$ and $Y$ by $\pi_X$ and $\pi_Y$.

Let $\cplx{K}$ be an object in $\dbc{X\times Y}$. The integral
functor defined by $\cplx{K}$ is the functor $\fmf{\cplx K}{X}{Y}
\colon \dmc{X} \to \dmc{Y}$ given by $$\fmf{\cplx
K}{X}{Y}(\cplx{F})=\bR \pi_{Y\ast}(\pi_X^\ast\cplx{F}\lotimes
\cplx{K})\, .$$ If the kernel $\cplx{K}\in \cdbc{X\times Y}$ is of
finite homological dimension over $X$, then the functor
$\fmf{\cplx K}{X}{Y}$ is defined over the whole $\catD(X)$ and
maps $\cdbc{X}$ to $\cdbc{Y}$.

If $Z$ is a third  proper scheme and $\cplx{L}$ is an object of
$\cdbc{Y\times Z}$, arguing exactly as in the smooth case, we
prove that there is an isomorphism of functors
\begin{equation}\label{e:convolution}
\fmf{\cplx{L}}{Y}{Z}\circ\fmf{\cplx{K}}{X}{Y}\simeq\fmf{\cplx{L}\ast\cplx{K}}{X}{Z} $$
where $$\cplx{L}\ast\cplx{K} = \bR \pi_{X,Z\ast}(\pi_{X,Y}^\ast
\cplx{K} \lotimes\pi_{Y,Z}^\ast\cplx{L})\,.
\end{equation}
If either $\cplx{K}$ or $\cplx{L}$ is of finite homological
dimension over $Y$, then $\cplx{L}\ast\cplx{K}$ in bounded.

\subsection{Adjoints}

We can describe nicely the adjoints to an integral functor when we work with Gorenstein schemes.
In this subsection $X$ and $Y$ are projective Gorenstein schemes.

\begin{prop} \label{p:AdjGo}  Let $\cplx{K}$ be an object in $\cdbc{X\times Y}$ of finite projective
dimension over $X$ and $Y$.
\begin{enumerate}
\item The
functor $\fmf{\dcplx{K}\otimes \pi_Y^\ast\omega_Y[n]}{Y}{X}
\colon \cdbc{Y}\to\cdbc{X}$ is a left adjoint  to the functor
$\fmf{\cplx{K}}{X}{Y}$.
\item The functor
$\fmf{\dcplx{K}\otimes
\pi_X^\ast\omega_X[m]}{Y}{X} \colon \cdbc{Y}\to\cdbc{X}$ is a
right adjoint  to the functor
$\fmf{\cplx{K}}{X}{Y}$.
\end{enumerate}
(here $m=\dim X$ and $n=\dim Y$)
\end{prop}
\begin{proof} We shall freely use \eqref{basechangedual} for the projections $\pi_X$ and $\pi_Y$.

(1) We first notice that one has
\begin{equation}
\label{e:formula} (\pi_X^\ast\cplx{F}\lotimes
\cplx{K})^{\scriptscriptstyle\vee}\simeq\bR\SHom_{\cO_{X\times
Y}}(\cplx{K},\pi_X^\ast\cplx{F}^{\scriptscriptstyle\vee})\simeq
\cplx{K}^{\scriptscriptstyle\vee}\lotimes\pi_X^\ast\cplx{F}^{\scriptscriptstyle\vee}\,,
\end{equation}
for $\cplx{F}$ in $\cdbc{X}$ by  \eqref{tens1} and Proposition
\ref{dualGore}. The latter applies because $\pi_X$ is a projective
morphism and $\cplx K$ is of finite projective dimension over $X$.
Now, if $\cplx G$ is an object of $\cdbc{Y}$ there is a chain of
isomorphisms
\begin{align*}
\Hom_{\catD(Y)}&(\cplx{G},\fmf{\cplx{K}}{X}{Y}(\cplx{F})) \simeq
\Hom_{\catD(X\times
Y)}(\pi_Y^\ast\cplx{G},\pi_X^\ast\cplx{F}\lotimes\cplx{K}) \\
&\simeq   \Hom_{\catD(X\times Y)}((\pi_X^\ast\cplx{F}\lotimes
\cplx{K})^{\scriptscriptstyle\vee},(\pi_Y^\ast\cplx{G})^{\scriptscriptstyle\vee})\\
&\simeq   \Hom_{\catD(X\times
Y)}(\pi_X^\ast\cplx{F}^{\scriptscriptstyle\vee}\lotimes
\cplx{K}^{\scriptscriptstyle\vee},\pi_Y^\ast\cplx{G}^{\scriptscriptstyle\vee})\\
&\simeq   \Hom_{\catD(X\times
Y)}(\pi_X^\ast\cplx{F}^{\scriptscriptstyle\vee},\bR\SHom_{\cO_{X\times
Y}}(\cplx{K}^{\scriptscriptstyle\vee},
\pi_Y^\ast\cplx{G}^{\scriptscriptstyle\vee}))\\
&\simeq \Hom_{\catD(X\times
Y)}(\pi_X^\ast\cplx{F}^{\scriptscriptstyle\vee},(\cplx{K}^{\scriptscriptstyle\vee} \lotimes
\pi_Y^\ast\cplx{G})^{\scriptscriptstyle\vee})\\
&\simeq
\Hom_{\catD(X)}(\cplx{F}^{\scriptscriptstyle\vee},\bR\pi_{X\ast}((\cplx{K}^{\scriptscriptstyle\vee} \lotimes
\pi_Y^\ast\cplx{G})^{\scriptscriptstyle\vee})) \\
&\simeq
\Hom_{\catD(X)}(\cplx{F}^{\scriptscriptstyle\vee}, \bR\SHom_{\cO_X}(\bR\pi_{X\ast}(\cplx{K}^{\scriptscriptstyle\vee}
\lotimes\pi_Y^\ast\cplx{G}\otimes\pi_Y^\ast\omega_Y[n]),\cO_{X}))\\
&\simeq \Hom_{\catD(X)}(\cplx{F}^{\scriptscriptstyle\vee}, [\fmf{\dcplx{K}\otimes
\pi_Y^\ast\omega_Y[n]}{Y}{X}(\cplx{G})]^\vee)  \\
&\simeq \Hom_{\catD(X)}(\fmf{\dcplx{K}\otimes
\pi_Y^\ast\omega_Y[n]}{Y}{X}(\cplx{G}),\cplx{F})\,.
\end{align*}
where the second follows from \eqref{homdual} which applies
because $\pi_X^\ast\cplx{F}\lotimes\cplx{K}$ is bounded, the third
is \eqref{e:formula}, the forth and the fifth are \eqref{tens1},
the seventh is relative duality and the ninth is again
\eqref{homdual}.

(2) The adjunction between the direct and inverse images and
relative duality proves that the functor
$$
H(\cplx{G})=\bR \pi_{X,\ast}( \bR\dSHom{\cO_{X\times Y}} (\cplx{K},
\pi_Y^!\cplx{G}))
$$
 satisfies
\begin{equation}  \label{e:adjunction}
\Hom_{\catD(Y)}(\fmf{\cplx{K}}{X}{Y}(\cplx{F}),\cplx{G}) \simeq
 \Hom_{\catD(X)}(\cplx{F},
 H(\cplx{G}))
\, .
\end{equation}
 Then  we conclude by Proposition \ref{dualGore} since $\pi_Y$ is a projective morphism.
\end{proof}

We shall need some basic results about adjoints and fully
faithfulness which we state without proof.
\begin{prop} \label{p:Adj} Let $\Phi \colon \cA\to\cB$ a functor and
$G\colon \cB\to \cA$ a left adjoint (resp.~$H\colon \cB\to
\cA$ a right adjoint). Then the following conditions are
equivalent:
\begin{enumerate}
\item $\Phi $ is fully faithful.
\item $G\circ \Phi $ is fully faithful (resp.~$H\circ \Phi $ is
fully faithful).
\item The counit morphism $G\circ \Phi \to \Id$ is an isomorphism
(resp.~ the unit morphism $\Id\to H\circ \Phi $ is an
isomorphism).
\end{enumerate}
Moreover, $\Phi $ is an equivalence if and only if $\Phi $ and $G$
(resp.~$\Phi $ and $H$) are fully faithful.
\qed
\end{prop}

\subsection{Strongly simple objects}

Let $X$ and $Y$ be proper Gorenstein schemes. In this situation,
the notion of strong simplicity is the following.

\begin{defn} An object  $\cplx{K}$ in $\cdbc{X\times Y}$ is \emph{strongly simple} over $X$ if it satisfies the
following conditions:
\begin{enumerate}
\item For every  closed point $x\in X$ there is a l.c.i.~zero cycle
$Z_x$ supported on $x$ such that
$$
\Hom^i_{\catD(Y)}(\fmf{\cplx{K}}{X}{Y}(\cO_{Z_{x_1}}),\fmf{\cplx{K}}{X}{Y}(\cO_{Z_{x_2}}))=0
$$
unless $x_1= x_2$ and $0\leq i\leq \dim X$.
\item$\Hom^0_{\catD(Y)}(\fmf{\cplx{K}}{X}{Y}(\sk{x}),
\fmf{\cplx{K}}{X}{Y}(\sk{x}))= k$ for every
%smooth
closed point $x\in X$.
\end{enumerate}
\label{1:strgsplcplxGo}
\end{defn}

The last condition can be written as $\Hom^0_{\catD(Y)}(\bL
j_{x}^\ast\cplx{K}, \bL j_{x}^\ast\cplx{K})= k$, because the
restriction $\bL j_x^\ast\cplx{K}$ of $\cplx{K}$ to the fibre $j_x\colon Y\simeq \{x\}\times Y\hookrightarrow X\times Y$ can
also be  computed  as $\fmf{\cplx{K}}{X}{Y}(\sk{x})$.

In order to fix some notation, for any zero-cycle $Z_x$ of $X$ and any scheme $S$, we shall denote by $j_{Z_x}$ the immersion $Z_x\times S \hookrightarrow X\times S$.

\begin{prop}\label{1:simpledual} Assume that $Y$ is projective, and let $\cplx K$ be a kernel in $\dbc{X\times Y}$ of finite projective dimension over $X$. If $\cplx K$ is strongly simple over $X$, its dual $\dcplx K$ is strongly simple over $X$ as well.
\end{prop}
\begin{proof} If $Z_x$ is a  l.c.i.~zero cycle supported on $x$, one has that $\fmf{\cplx{K}}{X}{Y}(\cO_{Z_x})
= p_{2\ast} \bL j_{Z_x}^\ast \cplx{K}$, with $p_2\colon Z_x\times
Y\to Y$ the second projection. Since $\omega_{p_2} \simeq
\cO_{Z_x\times Y}$ because $Z_x$ is zero dimensional and
Gorenstein, one obtains $$ \fmf{\cplx{K}}{X}{Y}(\cO_{Z_x})^\vee
\simeq p_{2\ast} (\bL j_{Z_x}^\ast \cplx{K})^{\vee}\,.$$ Moreover,
$(\bL j_{Z_x}^\ast \cplx{K})^{\vee}\simeq \bL j_{Z_x}^\ast
(\dcplx{K})$ by \eqref{basechangedual}
 since $j_{Z_x}$ is a regular closed immersion. Then,
$$
\fmf{\cplx{K}}{X}{Y}(\cO_{Z_x})^\vee \simeq \fmf{\dcplx{K}}{X}{Y}(\cO_{Z_x})\,.
$$
It follows that $\dcplx K$ satisfies condition (1) of Definition \ref{1:strgsplcplxGo}. To see that it also fulfils condition (1),
 we have to prove that $\fmf{\cplx{K}}{X}{Y}(\sk{x})^\vee \simeq \fmf{\cplx{K}^\vee}{X}{Y}(\sk{x})$, and this is equivalent to the
 base change formula $\bL j_x^\ast (\dcplx K) \simeq (\bL j_x^\ast \cplx K)^\vee$. Since we cannot longer use \eqref{basechangedual}
 because $j_x$ may fail to be of finite Tor-dimension, we proceed in a different way. 
 To see that the natural morphism $\bL j_x^\ast (\dcplx K) \to (\bL j_x^\ast \cplx K)^\vee$
 is an isomorphism, it suffices to check that
 $j_{x\ast}\bL j_x^\ast (\dcplx K) \simeq j_{x\ast}[(\bL j_x^\ast \cplx K)^\vee]$
 since $j_x$ is a closed embedding.
 On the one hand, we have
 $$
 j_{x\ast}\bL j_x^\ast (\dcplx K) \simeq 
 \dcplx K\lotimes j_{x\ast} \cO_Y \simeq
 \dcplx K\lotimes \pi_X^\ast \cO_x\,,
 $$
 whilst on the other hand,  
\begin{align*}
 j_{x\ast}[(\bL j_x^\ast \cplx K)^\vee] & =
 j_{x\ast} \bR\dSHom{\cO_{Y}}(\bL j_x^\ast \cplx K, \cO_Y) \simeq 
 \bR \dSHom{\cO_{X\times Y}}(\cplx K, j_{x\ast} \cO_Y) 
\\
&  \simeq 
    \bR\dSHom{\cO_{X\times Y}}(\cplx K, \pi_X^\ast \cO_x)\,.
 \end{align*}
We conclude by Proposition \ref{dualGore}.

% By Lemma \ref{l:beilinson} to prove
% that the natural morphism $\bL j_x^\ast (\dcplx K) \to (\bL j_x^\ast \cplx K)^\vee$  is an isomorphism, it suffices to see that
% $\bR \rho_{Y\ast} (\bL j_x^\ast (\dcplx K)(r)) \to \bR \rho_{Y\ast} ((\bL j_x^\ast \cplx K)^\vee(r))$ is an isomorphism for any $r$, where
% $\rho_Y\colon Y \to \Spec k$ is the projection on a point.  By base change and relative duality, one has
%$$
%\begin{aligned}
%\bR \rho_{Y\ast} (\bL j_x^\ast (\dcplx K)(r)) &\simeq  \bR \rho_{Y\ast} \bL j_x^\ast (\cplx K(-r)^\vee) \simeq
%\bL j_x^\ast \bR \pi_{X\ast} (\cplx K(-r)^\vee) \\
%& \simeq  \bL j_x^\ast  [(\bR \pi_{X\ast}(\cplx K(-r)\lotimes \pi_X^!\cO_X)]^\vee\,.
%\end{aligned}
%$$
%Analogously, one obtains
%$$
%\begin{aligned}
%\bR \rho_{Y\ast}((\bL j_x^\ast \cplx K)^\vee(r))& \simeq \bR  \rho_{Y\ast}((\bL j_x^\ast \cplx K(-r))^\vee ) \simeq [\bR \rho_{Y\ast}(\bL j_x^\ast \cplx K(-r)\lotimes \rho_Y^!\cO_x
%)]^\vee  \\
%& \simeq [\bR  \rho_{Y\ast} \bL j_x^\ast (\cplx K(-r)\lotimes \pi_X^!\cO_X )]^\vee \simeq
%[\bL j_x^\ast \bR \pi_{X\ast}(\cplx K(-r)\lotimes \pi_X^!\cO_X )]^\vee \,.
%\end{aligned}
%$$
%Now, since $\cplx K\lotimes \pi_X^!\cO_X$ has finite projective
%dimension over $X$,  the complex $\cplx M=\bR \pi_{X\ast} (\cplx
%K(-r)\lotimes \pi_X^!\cO_X)$ is of finite homological dimension,
%so $\bL j_x^\ast (\dcplx M) \simeq (\bL j_x^\ast \cplx M)^\vee$
%and we conclude.
\end{proof}

\begin{rema}
When $X$ and $Y$ are smooth, strong simplicity is usually defined by the
following conditions (see \cite{BBHJ}):
\begin{enumerate}
\item $\Hom^i_{\catD(Y)}(\bL j_{x_1}^\ast\cplx{K}, \bL
j_{x_2}^\ast\cplx{K})=0$ unless $x_1= x_2$ and $0\leq i\leq \dim
X$; \item $\Hom^0_{\catD(Y)}(\bL j_{x}^\ast\cplx{K}, \bL
j_{x}^\ast\cplx{K})= k$ for every closed point $x$.
\end{enumerate}
Since our definition is weaker, Theorem \ref{1:ffcritGo} improves
Bondal and Orlov's result \cite[Thm.~1.1]{BO95}.
\end{rema}

We now give the criterion for an integral functor between derived
categories of Gorenstein proper schemes to be fully faithful.

\begin{thm}\label{1:ffcritGo}  Let $X$ and $Y$ be projective Gorenstein schemes over an algebraically closed field of characteristic zero,
and let $\cplx{K}$ be an object in $\cdbc{X\times Y}$ of finite
projective dimension over $X$ and over $Y$. Assume also that $X$
is
%irreducible
integral. Then the functor $\fmf{\cplx{K}}{X}{Y}\colon
\cdbc{X}\to \cdbc{Y}$ is fully faithful if and only if the kernel
$\cplx{K}$ is strongly simple over $X$.
\end{thm}
\begin{proof}  If the functor $\fmf{\cplx{K}}{X}{Y}$ is fully faithful, then $\cplx{K}$ is
strongly simple over $X$.

Let us prove the converse. Before starting, we fix some notation:
we denote by $\pi_i$ the projections of $X\times X$ onto its
factors and by $U$ the smooth locus of $X$, which is not empty
because $X$ is integral. We also denote $m=\dim X$, $n=\dim Y$ and
$\Phi =\fmf{\cplx{K}}{X}{Y}$.

By Proposition \ref{p:AdjGo}, $\Phi $ has a left adjoint
$G=\fmf{\dcplx{K}\otimes \pi_Y^\ast\omega_Y[n]}{Y}{X}$  and a
right adjoint $H=\fmf{\dcplx{K}\otimes
\pi_X^\ast\omega_X[m]}{Y}{X}$. By Proposition \ref{p:Adj} it
suffices to show that $G\circ \Phi $ is fully faithful. We know
that $H\circ \Phi  \simeq \fmf{\cplx{M}}{X}{X}$, and $G\circ \Phi
\simeq \fmf{\widetilde{\cplx{M}}}{X}{X}$, with ${\cplx M}$ and
$\widetilde{\cplx{M}}$ given by \eqref{e:convolution}. Notice that
since $\cplx{K}$ is of finite projective dimension over $X$ and
$Y$, $\cplx{M}$ and $\widetilde{\cplx{M}}$ are bounded.

The strategy of the proof is as follows: we first show that both
$\cplx{M}$ and $\widetilde{\cplx M}$ are single sheaves supported
topologically on the image $\Delta$ of the diagonal morphism
$\delta \colon X\hookrightarrow X\times X$; then we prove that
$\widetilde{\cplx M}$ is actually schematically supported on the
diagonal, that is, $\widetilde{\cplx{M}} = \delta_\ast \cN$ for a
coherent sheaf $\cN$ on $X$ and finally that $\cN$ is a line
bundle; this will imply that $\fmf{\widetilde{\cplx{M}}}{X}{X}$ is the twist by $ \cN$ which is an equivalence of categories, in particular fully faithful.

\medskip\noindent
a) \emph{$\cplx{M}$ and $\widetilde{\cplx M}$ are single sheaves topologically supported on the diagonal.}

\medskip

Let us fix a closed point $(x_1,x_2)\in X\times X$ and consider the
l.c.i.~zero cycles $Z_{x_1}$ and $Z_{x_2}$ of the first condition of
the definition of strongly simple object.  One has
$$
\Hom^i_{\catD(X)}(\cO_{Z_{x_1}},
\fmf{\cplx{M}}{X}{X}(\cO_{Z_{x_2}})) \simeq \Hom^i_{\catD(Y)}(\Phi(\cO_{Z_{x_1}}),\Phi(\cO_{Z_{x_2}}))\,,
$$
which is zero unless $x_1=x_2$ and $0\le i\le m$ because $\cplx K$
is strongly simple.  Applying Proposition \ref{1:support2} to  the
immersion $\{x_2\}\hookrightarrow X$ we have that
$\fmf{\cplx{M}}{X}{X}(\cO_{Z_{x_2}})$ reduces to a coherent sheaf
topologically supported at $x_2$. Since
$\fmf{\cplx{M}}{X}{X}(\cO_{Z_{x_2}})\simeq p_{2\ast}\bL
j_{Z_{x_2}}^* \cplx{M} $, where $p_2\colon Z_{x_2}\times X\to X$
is the second projection, the complex $\bL j_{Z_{x_2}}^* \cplx{M}$
is isomorphic to a single coherent sheaf $\cF$ topologically
supported at $(x_2,x_2)$. If we denote by $i_{Z_{x_1}}\colon
Z_{x_2}\times Z_{x_1}\hookrightarrow Z_{x_2}\times X$ and
$j_{Z_{x_2}\times Z_{x_1}}\colon Z_{x_2}\times
Z_{x_1}\hookrightarrow X\times X$ the natural immersions, we have
$$
\bL j_{Z_{x_2}\times Z_{x_1}}^*\cplx{M}\simeq \bL i_{Z_{x_1}}^* \bL j_{Z_{x_2}}^* \cplx{M} \simeq \bL i_{Z_{x_1}}^* \cF\,.
$$
Thus, $L_p j_{Z_{x_2}\times Z_{x_1}}^*\cplx{M}=0$ unless
$x_1=x_2$ and $0\le p\le m$. Applying now Proposition
\ref{1:support} to $\delta$, we obtain that $\cplx{M}$ reduces to a coherent sheaf $\cM$ supported topologically on the diagonal as claimed.

For $\widetilde{\cplx M}$, we proceed as follows. We have
\begin{align*}
\calH^i(\bL j_{Z_{x_2}}^\ast \fmf{\widetilde{\cplx{M}}}{X}{X}(\cO_{Z_{x_1}})^\vee) & \simeq
\Hom^i_{\catD(Z_{x_2})}(\bL
j_{Z_{x_2}}^\ast \fmf{\widetilde{\cplx{M}}}{X}{X}(\cO_{Z_{x_1}}),
\cO_{Z_{x_2}}) \\
& \simeq \Hom^i_{\catD(X)}(\fmf{\widetilde{\cplx{M}}}{X}{X}
(\cO_{Z_{x_1}}), \cO_{Z_{x_2}})
 \simeq
\Hom^i_{\catD(Y)}(\Phi(\cO_{Z_{x_1}}),\Phi(\cO_{Z_{x_2}})),
\end{align*}
which is zero unless $x_1=x_2$ and $0\le i\le m$ because $\cplx K$ is strongly simple.
Since $Z_{x_2}$ is a zero dimensional Gorenstein scheme,  \eqref{e:zeroddual} implies that
$L_i j_{Z_{x_2}}^\ast \fmf{\widetilde{\cplx{M}}}{X}{X}(\cO_{Z_{x_1}})=0$ again unless $x_1=x_2$ and $0\le i\le m$.  By Proposition \ref{1:support} for the immersion $\{x_1\}\hookrightarrow X$, one has that $
\fmf{\widetilde{\cplx{M}}}{X}{X}(\cO_{Z_{x_1}})$ is a sheaf supported
topologically at $x_1$. Now, a similar argument to the one used for $\cplx{M}$ proves that $\widetilde{\cplx{M}}$ reduces to a coherent sheaf
$\widetilde\cM$ supported topologically on the diagonal.

\medskip\noindent
b) \emph{$\widetilde{\cM}$ is schematically supported on the
diagonal, that is, $\widetilde{\cM} = \delta_\ast \cN$ for a
coherent sheaf $\cN$ on $X$; moreover $\cN$ is a line bundle.}

\medskip
It might happen that the schematic support is an infinitesimal
neighborhood of the diagonal; we shall see that this is not the
case. Let us denote by $\bar\delta\colon W \hookrightarrow X\times
X$ the schematic support of $\widetilde\cM$ so that $\widetilde\cM
= \bar \delta_\ast \cN$ for a coherent sheaf $\cN$ on $W$. The
diagonal embedding $\delta$ factors through a closed immersion
$\tau\colon X\hookrightarrow W$ which topologically is a
homeomorphism.

\medskip\noindent

b1) \emph{$\pi_{2\ast}(\widetilde\cM)$ is locally free.}

\medskip

 To see this, we shall prove that $\Hom_{\catD(X)}^1(\pi_{2\ast}(\widetilde\cM),\cO_x)=0$ for every closed point $x\in X$.
Since $\widetilde\cM$ is topologically supported on the diagonal,
we have that $\pi_{2\ast}(\widetilde\cM)\simeq \bR
\pi_{2\ast}(\widetilde\cM)\simeq \fmf{\widetilde{\cplx
M}}{X}{X}(\cO_X)$. We have
$$
\Hom^1_{\catD(X)}(\pi_{2\ast}(\widetilde{\cM}),\cO_x)\simeq
\Hom^1_{\catD(X)}(\fmf{\widetilde{\cplx M}}{X}{X}(\cO_X),\cO_x)
 \simeq \Hom^1_{\catD(X)}(\cO_X, \fmf{\cplx M}{X}{X}(\cO_x))
 $$
for every closed point $x\in X$, because $H\circ \Phi$ is a
right adjoint to $G\circ \Phi$. Since $\fmf{\cplx M}{X}{X}(\cO_x)\simeq \bL j_x^\ast \cM$ has only
negative cohomology sheaves and all of them are supported at $x$, one has that  $\Hom_{\catD(X)}^1(\pi_{2\ast}(\widetilde\cM),\cO_x)=0$
and $\pi_{2\ast}(\widetilde\cM)$ is locally free.

\medskip\noindent

b2) \emph{$\pi_{1\ast}(\widetilde\cM)$ is a line bundle on the smooth locus $U$ of $X$.}

We know that $\fmf{\cM}{X}{X}(\cO_{Z_{x_2}})$ reduces to a single sheaf supported at $x_2$. Then, for every point $x_1\in U$ one has
\begin{align*}
\calH^i(\bL j_{Z_{x_2}}^\ast \fmf{\widetilde{\cplx{M}}}{X}{X}(\cO_{x_1})^\vee) & \simeq
\Hom^i_{\catD(Z_{x_2})}(\bL
j_{Z_{x_2}}^\ast \fmf{\widetilde{\cplx{M}}}{X}{X}(\cO_{x_1}),
\cO_{Z_{x_2}}) \\
& \simeq \Hom^i_{\catD(X)}(\fmf{\widetilde{\cplx{M}}}{X}{X}
(\cO_{x_1}), \cO_{Z_{x_2}})
 \simeq
\Hom^i_{\catD(X)}(\cO_{x_1},\fmf{\cM}{X}{X}(\cO_{Z_{x_2}}))
\end{align*}
which is zero unless $x_2=x_1$ and $0\le i\le m$ because $x_1$ is a smooth point. Since $Z_{x_2}$ is a zero dimensional Gorenstein scheme, \eqref{e:zeroddual} implies that whenever $x_1$ is a smooth point, then $L_i j_{Z_{x_2}}^\ast \fmf{\widetilde{\cplx{M}}}{X}{X}(\cO_{x_1})=0$ unless $x_2=x_1$ and $0\le i\le m$. By Proposition \ref{1:support}, $\fmf{\widetilde{\cplx{M}}}{X}{X}(\cO_{x_1})$ reduces to a single sheaf supported at $x_1$. In particular $\bL j_{x}^\ast (\widetilde \cM)\simeq j_{x}^\ast (\widetilde \cM)$ for every smooth point $x$, and thus the restriction of $\widetilde\cM$ to $U\times X$ is flat over $U$.
Moreover, for every point $x\in U$, we have that
$$
\Hom_X(j_x^\ast\widetilde \cM, \cO_x)\simeq \Hom^0_{\catD (X)}(\Phi(\cO_x), \Phi(\cO_x))\simeq k\,.
$$
By \cite[Lemmas 5.2 and 5.3]{Bri99} there is a point $x_0$ in $U$ such that the Kodaira-Spencer map for the family $\rest{\widetilde\cM}{U\times X}$ is injective at $x_0$. We now proceed as in the proof of \cite[Thm.~5.1]{Bri99}: the morphism $\Hom_{\catD(X)}^1(\cO_{x_0},\cO_{x_0})\to \Hom_{\catD(X)}^1(G\circ \Phi (\cO_{x_0}), G\circ \Phi (\cO_{x_0}))$ is injective so that the morphism $\Hom_{\catD(X)}^1(\cO_{x_0},\cO_{x_0})\to \Hom_{\catD(Y)}^1(\Phi (\cO_{x_0}),\Phi (\cO_{x_0}))$ is injective as well and then the counit morphism $j_{x_0}^\ast\widetilde\cM\simeq G\circ \Phi (\cO_{x_0}) \to \cO_{x_0}$ is an isomorphism. Thus, the rank of $\pi_{1\ast}(\widetilde\cM)$ at the point $x_0$ is one, and then it is one everywhere in $U$.

\medskip\noindent

b3) \emph{$\rest{\tau}{U}\colon U\hookrightarrow W_U=W\cap (U\times X)$ is an isomorphism and $\cN'=\rest{\cN}{U}$ is a line bundle.}

\medskip
We proceed locally. We then write $U=\Spec A$, $W_U=\Spec B$ so that
$\tau$ corresponds to a surjective ring morphism $B\to A\to 0$ and the
projection $q_1=\rest{\pi_1}{U}\colon W_U\to U$ to an immersion
$A\hookrightarrow B$. Now $\cN'$ is a $B$-module which is
isomorphic to $A$ as an $A$-module, $\cN'\simeq e\cdot A$,  because
$ q_{1\ast}(\cN')=\rest{\pi_{1\ast}(\widetilde\cM)}{U}$ is a line bundle. It follows
that $\cN'$ is also generated by $e$  as a $B$-module. The kernel
of $B\to \cN'\simeq e\cdot B\to 0$ is the annihilator of $\cN'$ and
then it is zero by the very definition of $W$. It follows that
$B\simeq A$ as an $A$-module and then the morphism $B\to A\to 0$
is an isomorphism. Hence, $W_U\simeq U$, $q_1$ is the identity
map, and $\cN'\simeq q_{1\ast}(\cN')$ is a line bundle.

\medskip\noindent

b4) \emph{$\tau\colon X\hookrightarrow W$ is an isomorphism and $\cN$ is a line bundle.}

\medskip
Since $U\simeq W_U$, $\rest{\pi_{2\ast}\widetilde \cM}{U}\simeq \rest{\cN}{U}\simeq  \rest {\pi_{1\ast}\widetilde \cM}{U}$, which is a line bundle on $U$. Then, the locally free sheaf  $\pi_{2\ast}\widetilde \cM$ has to be a line bundle. Then the same argument used in b3)  proves the remaining statement.
\end{proof}

\begin{cor}  \label{c:strgsplcplxGo}  An object  $\cplx{K}$ in $\cdbc{X\times Y}$ satisfying the
conditions of Theorem \ref{1:ffcritGo} is strongly simple over $X$ if and only if
\begin{enumerate}
\item
$\Hom^i_{\catD(Y)}(\fmf{\cplx{K}}{X}{Y}(\cO_{Z_{x_1}}),\fmf{\cplx{K}}{X}{Y}(\cO_{Z_{x_2}}))=0$
for any pair  $Z_{x_1}$ and $Z_{x_2}$ of l.c.i.~zero cycles
(supported on $x_1$, $x_2$ respectively) unless $x_1= x_2$ and
$0\leq i\leq \dim X$;
\item$\Hom^0_{\catD(Y)}(\fmf{\cplx{K}}{X}{Y}(\sk{x}),
\fmf{\cplx{K}}{X}{Y}(\sk{x}))= k$ for every point $x\in X$.
\end{enumerate}
\end{cor}

From Propositions \ref{1:simpledual} and \ref{p:Adj} and Corollary \ref{fpddual}, we
obtain:

\begin{cor} \label{c:ffcritGo} Let $X$ and $Y$ be projective  integral
%irreducible
Gorenstein
schemes over an algebraically closed field of characteristic zero, and let $\cplx{K}$ be an object in $\cdbc{X\times Y}$ of
finite projective dimension over both factors. The integral
functor $\fmf{\cplx{K}}{X}{Y}$ is an equivalence if and only if
$\cplx{K}$ is strongly simple over both factors.
\end{cor}

\begin{rem} \label{r:charp} Theorem \ref{1:ffcritGo} is false in positive characteristic even in the smooth case. Let $X$ be a smooth projective scheme of dimension $m$ over a field $k$ of characteristic $p>0$,
%We have a commutative diagram
%$$
%\xymatrix{
%X \ar[r]^{F_X} \ar[d] & X \ar[d] \\
%\Spec k \ar[r]^{F_k} & \Spec k
%}
%$$
%where $F_k$ and $F_X$ are the corresponding Frobenius morphisms. We denote by $X^{(p)}$ the fibre product of the morphisms $F_k\colon \Spec k \to \Spec k$ and $X \to \Spec k$. The morphism $F_X$ defines a morphism of $k$-schemes
and
$F\colon X \to X^{(p)}$
the \emph{relative Frobenius morphism} \cite[3.1]{Illu96}, which is topologically a homeomorphism.  Let $\Gamma\hookrightarrow X\times X^{(p)}$ be the graph of $F$, whose associated integral functor is the direct image $F_\ast \colon \cdbc{X} \to \cdbc{X^{(p)}}$. Since $F_\ast(\cO_x)\simeq \cO_{F(x)}$, one easily sees that $\Gamma$ is strongly simple over $X$. However, $F_\ast(\cO_{X})$ is a locally free $\cO_{X^{(p)}}$-module of rank $p^m$ \cite[3.2]{Illu96}, so that $\Hom^0_{\catD(X^{(p)})}(F_\ast(\cO_X), \cO_{F(x)})\simeq k^{p^m}$ whereas $\Hom^0_{\catD(X)}(\cO_{X},\cO_x)\simeq k$; thus $F_\ast$ is not fully faithful.
\end{rem}
\subsection{A criterion for equivalence}

The usual Bridgeland criterion \cite[Thm.~5.1]{Bri99} that
characterises when an integral functor over the derived category
of a smooth variety is an equivalence (or a Fourier-Mukai functor)
also works in the Gorenstein case. The original proof is based on
the fact that if $X$ is smooth, the skyscraper sheaves $\cO_x$
form a spanning class for the derived category $\cdbc{X}$
\cite{Bri99}. This is also true for Gorenstein varieties.
%the
%only thing that  one always has is that
%\begin{equation}
%\Hom^i(\cplx{E},\cO_x)=0 \text{ for every $i$ and every $x\in X$} \implies \cplx{E}=0\,.
%\label{e:spann}
%\end{equation}
Moreover in this case there is a more natural spanning class (see
\cite{Chen02} for a similar statement), that allows to give a
similar criterion.
\begin{lem} \label{l:spanning} If $X$ is a Gorenstein scheme, then
the following sets are spanning classes for $\cdbc{X}$:
\begin{enumerate}
\item $\Omega_1=\{\cO_x\}$ for all closed points $x\in X$.
\item $\Omega_2=\{ \cO_{Z_x}\}$ for all closed points $x\in X$ and all l.c.i.~zero cycles $Z_x$ supported on $x$.
\end{enumerate}
\end{lem}
\begin{proof} (1) Arguing as in \cite[Lemma 2.2]{Bri99}, one proves
that if $\Hom^i(\cplx{E},\cO_x)=0$ for every $i$ and every $x\in
X$, then  $\cplx{E}=0$.
Suppose now that $\Hom^i(\cO_x,\cplx{E})=0$ for every $i$ and
every $x\in X$. By \eqref{homdual},
$\Hom^i(\cO_x,\cplx{E})\simeq\Hom^i(\dcplx{E},\cO_x^\vee)$ and
since $\cO_x^\vee\simeq \cO_x[-m]$ where $m=\dim X$ because $X$ is
Gorenstein, we have that $\Hom^{i-m}(\dcplx{E},\cO_x)=0$ for every
$i$ and every $x\in X$. Then $\dcplx{E}=0$ and from
\eqref{ddual}, one concludes that $\cplx{E}=0$.

(2) By Proposition \ref{1:support2} with $Y=\emptyset$, if
$\Hom^{i}_{\catD (X)}(\cO_{Z_x},\cplx E)= 0$ for every $i$ and
every $Z_x$, then $\cplx{E}=0$.
On the other hand, since $\cO_{Z_x}$ is of finite homological dimension, Serre duality can be applied to get
an isomorphism
$$
\Hom^i(\cplx E,\cO_{Z_x})^\ast\simeq \Hom^{i}(\cplx
E,\cO_{Z_x}\otimes \omega_X)^\ast\simeq \Hom^{m-i}(\cO_{Z_x},
\cplx E)
$$
where $m=\dim X$. By the first part, if $\cplx{E}$ is a non-zero
object in $\cdbc{X}$ the second member is non-zero for some $i$
and we finish.
\end{proof}
\begin{thm} \label{t:equivalence0}
Let $X$ and $Y$ be projective Gorenstein schemes over an algebraically closed field of characteristic zero. Assume also that
$X$ is
integral
%irreducible
 and $Y$ is connected. If $\cplx{K}$ is an
object in $\cdbc{X\times Y}$ of finite projective dimension over
both $X$ and $Y$, then the functor $\fmf{\cplx{K}}{X}{Y}\colon
\cdbc{X}\to \cdbc{Y}$ is an equivalence of categories if and only
if one has
\begin{enumerate}
\item $\cplx{K}$ is strongly simple over $X$. \item For every
closed point $x\in X$, $\fmf{\cplx K}{X}{Y}(\cO_x) \simeq
\fmf{\cplx K}{X}{Y}(\cO_x)\otimes \omega_Y$.
\end{enumerate}
\end{thm}
\begin{proof} By Proposition \ref{p:AdjGo}, the functor $H=\fmf{\dcplx{K}\otimes \pi_X^\ast\omega_X[m]}{Y}{X}$
is a right adjoint to $\fmf{\cplx K}{X}{Y}$ while
$G=\fmf{\dcplx{K}\otimes \pi_Y^\ast\omega_Y[n]}{Y}{X}$ is a left
adjoint to it. If $\fmf{\cplx K}{X}{Y}$ is an equivalence,  there
is an isomorphism of functors  $H\simeq G$ and then the left
adjoints are also isomorphic, that is $\fmf{\cplx K}{X}{Y}\simeq
\fmf{\cplx K\otimes \pi_Y^\ast\omega_Y^{-1} \otimes
\pi_X^\ast\omega_X}{X}{Y}$. Applying this to $\cO_x$ we get
$\fmf{\cplx K}{X}{Y}(\cO_x)\otimes \omega_Y\simeq \fmf{\cplx
K}{X}{Y}(\cO_x)$.

For the converse, notice first that  the derived category
$\dbc{Y}$ is indecomposable  because $Y$ is connected
\cite[Ex.~3.2]{Bri99}. Then we have to prove that for any object
$\cplx E$ in $\cdbc{Y}$ the condition $H(\cplx E)=0$ implies that
$G(\cplx E)=0$ \cite[Thm.~3.3]{Bri99}. Since for every object
$\cplx E$ in $\cdbc{Y}$ one has a functorial isomorphism
\begin{equation}
G(\cplx{E})\simeq H(\cplx{E}\otimes \omega_Y[n])\otimes \omega_X^{-1}[-m]\,,
\label{e:relacionadj}
\end{equation}
it is enough to prove that $H(\cplx E\otimes\omega_Y[n])=0$. We have
$$
\begin{aligned}
\Hom^i(\cO_x,H(\cplx{E}\otimes\omega_Y[n]))&\simeq
\Hom^i(\fmf{\cplx K}{X}{Y}(\cO_x), \cplx{E}\otimes\omega_Y[n])\\
&\simeq \Hom^{n+i}(\fmf{\cplx K}{X}{Y}(\cO_x), \cplx{E})\\
&\simeq \Hom^{n+i}(\cO_x,H(\cplx{E}))=0
\end{aligned}
$$
and one concludes by Lemma \ref{l:spanning}.
\end{proof}
Using now the second part of Lemma \ref{l:spanning}, we prove
analogously the following:

\begin{thm} \label{t:equivalence}
Let $X$ and $Y$ be projective Gorenstein schemes over an
algebraically closed field of characteristic zero. Assume also
that $X$ is integral
%irreducible
 and $Y$ is connected. If $\cplx{K}$ is an
object in $\cdbc{X\times Y}$ of finite projective dimension over
both $X$ and $Y$, then the functor $\fmf{\cplx{K}}{X}{Y}\colon
\cdbc{X}\to \cdbc{Y}$ is an equivalence of categories if and only
if one has
\begin{enumerate}
\item $\cplx{K}$ is strongly simple over $X$. \item For every
closed point $x\in X$ there is a l.c.i.~cycle $Z_x$ such that
$\fmf{\cplx K}{X}{Y}(\cO_{Z_x}) \simeq \fmf{\cplx
K}{X}{Y}(\cO_{Z_x})\otimes \omega_Y$. \end{enumerate}
\end{thm}
\begin{rem} The second condition in the above lemma can be also written in either the form
$p_{2\ast}(\bL j_{Z_x}^\ast\cplx K)\simeq p_{2\ast} (\bL
j_{Z_x}^\ast\cplx K) \otimes  \omega_Y$ or the form $\bL
j_{Z_x}^\ast\cplx K\simeq  \bL j_{Z_x}^\ast\cplx K \otimes
p_2^\ast \omega_Y$, where $p_2\colon Z_x\times Y \to Y$ is the
projection.
\end{rem}

\subsection{Geometric applications of Fourier-Mukai functors}

As in the smooth case, the existence of a Fourier-Mukai functor between the derived categories of two
Gorenstein schemes has important geometrical consequences. In the following proposition, we list some of them.

\begin{prop} \label{p:conseq}
Let $X$ and $Y$ be projective Gorenstein schemes and let
$\cplx{K}$ be an object in $\cdbc{X\times Y}$ of finite projective
dimension over both $X$ and $Y$. If the integral functor
$\fmf{\cplx{K}}{X}{Y}\colon \cdbc{X}\to \cdbc{Y}$ is a
Fourier-Mukai functor, the following statements hold:
\begin{enumerate} \item The right and the left adjoints to
$\fmf{\cplx{K}}{X}{Y}$ are functorially isomorphic
$$
\fmf{\dcplx{K}\otimes \pi_X^\ast\omega_X[m]}{Y}{X}\simeq \fmf{\dcplx{K}\otimes
\pi_Y^\ast\omega_Y[n]}{Y}{X}
$$
and they are quasi-inverses
to $\fmf{\cplx{K}}{X}{Y}$. \item $X$ and $Y$ have the same
dimension, that is, $m=n$.
\item $\omega_X^r$ is trivial for an integer $r$ if and only
if $\omega_Y^r$ is trivial. Particularly, $\omega_X$ is trivial if and only
if $\omega_Y$ is trivial. In this case, the functor
$\fmf{\dcplx{K}}{Y}{X}$ is a quasi-inverse to
$\fmf{\cplx{K}}{X}{Y}$.
\end{enumerate}
\end{prop}

\begin{proof} (1) Since $\fmf{\cplx{K}}{X}{Y}$ is an equivalence, its quasi-inverse is a right and a left
adjoint. The statement follows from Proposition \ref{p:AdjGo}
using the uniqueness of the adjoints.

(2) Applying the above isomorphism to $\cO_{Z_y}$ where $Z_y$ is a
l.c.i.~zero cycle supported on $y$,  one obtains $p_{X\ast}(\bL j_{Z_y}^\ast\dcplx
K)[n]\simeq p_{X\ast} (\bL j_{Z_y}^\ast\dcplx K) \otimes
\omega_X[m]$ where $p_X\colon X\times Z_y\to X$ is the projection.
Since the two functors are equivalences, these are non-zero
objects in $\cdbc{X}$. Let $q_0$ be the minimum (resp.~maximum) of
the $q$'s with $\calH^q(p_{X\ast}(\bL j_{Z_y}^\ast\dcplx K))\neq
0$. Since $\calH^{q_0}( p_{X\ast}(\bL j_{Z_y}^\ast\dcplx K))\simeq
\calH^{q_0+m-n}(p_{X\ast}(\bL j_{Z_y}^\ast\dcplx K))\otimes
\omega_X$ one has $\calH^{q_0+m-n}(p_{X\ast}(\bL
j_{Z_y}^\ast\dcplx K))\otimes \omega_X\neq 0$ which contradicts
the minimality (resp.~maximality) if $m-n<0$ (resp.~$>0$). Thus,
$n=m$.

(3) If we denote
by $H$  the right adjoint to
$\fmf{\cplx{K}}{X}{Y}$, thanks to (1) and
\eqref{e:relacionadj} we have
that $H(\cplx{E})\otimes \omega_X^r\simeq H(\cplx{E}\otimes
\omega_Y^r)$ for every $\cplx E\in \cdbc{Y}$ and every integer $r$.  If $\omega_X^r\simeq\cO_X$, taking
$\cplx E=\cO_Y$ we have $H(\cO_Y)\simeq H(\omega_Y^r)$ and applying the functor $\fmf{\cplx{K}}{X}{Y}$ to
this isomorphism we get $\omega_Y^r\simeq \cO_Y$. The converse is similar.
\end{proof}

\section{Relative Fourier-Mukai transforms for Gorenstein morphisms}\label{s:relative}
\subsection{Generalities and base change properties}
Let $S$ be a
%reduced
scheme and let $p\colon X\to S$
and $q\colon Y\to S$ be proper morphisms. We denote by $\pi_X$ and
$\pi_Y$ the projections of the fibre product $X\times_SY$ onto its
factors and by $\rho=p\circ\pi_X=q\circ \pi_Y$ the projection of
$X\times_SY$ onto the base scheme $S$ so that we have the
following cartesian diagram
$$\xymatrix{ & X\times_S Y \ar[ld]_{\pi_{X}}\ar[rd]^{\pi_{Y}} \ar[dd]^{\rho}&  \\
X \ar[rd]^{p} & & Y\ar[ld]_{q} \\
& S & }$$
Let $\cplx{K}$ be an object in $\dbc{X\times_SY}$. The relative
integral functor defined by $\cplx{K}$ is the functor $\fmf{\cplx
K}{X}{Y} \colon \dmc{X} \to \dmc{Y}$ given by $$\fmf{\cplx
K}{X}{Y}(\cplx{F})=\bR \pi_{Y\ast}(\bL\pi_X^\ast\cplx{F}\lotimes
\cplx{K})\, .$$ We shall denote this functor by $\Phi$ from now on.

Let $s\in S$ be a closed point. Let us denote $X_s=p^{-1}(s)$,
$Y_s=q^{-1}(s)$, and $\Phi_s\colon \dmc{X_s} \to \dmc{Y_s}$ the
integral functor defined by $\cplx{K}_s=\bL j_s^\ast \cplx{K}$,
with $j_s\colon X_s\times Y_s\hookrightarrow X\times_S Y$ the
natural embedding.

When the kernel $\cplx{K}\in \cdbc{X\times_SY}$ is of finite
homological dimension over $X$, the functor $\Phi$ is defined over
the whole $\catD(X)$ and it maps $\cdbc{X}$ into $\cdbc{Y}$. If
moreover $q\colon Y\to S$ is flat, then $\cplx{K}_s$ is of finite
homological dimension over $X_s$ for any $s\in S$.

If $p\colon X\to S$ and $q\colon Y\to S$ are flat morphisms, from the base-change
formula we obtain that
\begin{equation}\label{e:basechange}\bL
j_s^\ast\Phi (\cplx{F})\simeq\Phi_s(\bL j_s^\ast\cplx{F})
\end{equation}
for every $\cplx{F}\in \catD(X )$, where $j_s\colon
X_s\hookrightarrow X$ and $j_s\colon Y_s\hookrightarrow Y$ are the
natural embeddings.  In this situation, base change formula also gives that
\begin{equation}\label{e:directimage}
j_{s\ast}\Phi_s(\cplx{G})\simeq\Phi ( j_{s\ast}\cplx{G})
\end{equation}
for every $\cplx{G}\in \catD(X_s)$.

Proposition \ref{dualGore} allows us to obtain the following
result.

\begin{lem} \label{relatrightadj}Let $p\colon X\to S$ and $q\colon Y\to S$ be locally projective
Gorenstein morphisms, and let $\cplx{K}$ be an object in
$\cdbc{X\times_SY}$ of finite projective dimension over both $X$ and $Y$. Then the functor
$$
H=\fmf{\dcplx{K}\lotimes \pi_X^\ast\omega_{X/S}[m]}{Y}{X}
\colon \cdbc{Y}\to\cdbc{X}$$ is a right adjoint  to the functor
$\fmf{\cplx{K}}{X}{Y}$.
\end{lem}

\subsection{Criteria for fully faithfulness and equivalence in the relative setting}

In this subsection we work over an algebraically closed field of characteristic zero.

In the relative situation the notion of strongly simple object is the following.

\begin{defn} Let $p\colon X\to S$ and $q\colon Y\to S$ be proper
Gorenstein morphisms. An object $\cplx{K}\in \cdbc{X\times_SY}$
is \emph{relatively strongly simple} over $X$ if $\cplx{K}_s$ is bounded and strongly simple over $X_s$  for every closed point $s\in S$.
\end{defn}

\begin{lem} \label{l:fibres} Let $Z\to S$ be a proper morphism and $\cplx E$ be an object of $\cdbc{Z}$
such that $\bL j_s^\ast\cplx E=0$ in $\cdbc{Z_s}$ for every closed point $s$ in $S$, where
$j_s\colon Z_s\hookrightarrow Z$ is the immersion of the fibre. Then $\cplx E=0$.
\end{lem}
\begin{proof} For every closed point $s$ in $S$ there is a spectral sequence $E_2^{-p,q}=L_p j_s^\ast
\calH^q(\cplx E)$ converging to $E_\infty^{p+q}=\calH^{p+q}(\bL
j_s^\ast\cplx E)=0$. Assume that $\cplx E\neq 0$ and let $q_0$ be
the maximum of the integers $q$ such that $\calH^q(\cplx E)\neq
0$. If $s$ is a point in the image of the support of
$\calH^{q_0}(\cplx E)$, one has that $j_s^\ast \calH^{q_0}(\cplx
E)\neq 0$ and every non-zero element in $E_2^{0,q_0}=j_s^\ast
\calH^{q_0}(\cplx E)$ survives to infinity. Then
$E_\infty^{q_0}\neq 0$ and this is impossible.
\end{proof}
\begin{thm}  \label{t:relative}
Let $p\colon X\to S$ and $q\colon Y\to S$ be locally projective
Gorenstein morphisms. Let $\cplx{K}$ be an object in
$\cdbc{X\times_SY}$ of finite projective dimension over both $X$ and $Y$. The relative integral
functor $\Phi=\fmf{\cplx{K}}{X}{Y}\colon \cdbc{X}\to \cdbc{Y}$ is
fully faithful (resp.~an equivalence) if and only if $\Phi_s\colon
\cdbc{X_s}\to \cdbc{Y_s}$ is fully faithful (resp.~an equivalence)
for every closed point $s\in S$.
\end{thm}

\begin{proof} By  Proposition \ref{p:Adj}, if $\Phi$ is fully faithful the unit morphism
$$
\Id \to H\circ \Phi
$$
is an isomorphism (where $H$ is the right adjoint given at Lemma
\ref{relatrightadj}). Then, given a closed point $s\in S$ and
$\cplx G\in \cdbc{X_s}$, one has an isomorphism $j_{s\ast}\cplx G\to
(H\circ \Phi)(j_{s\ast}\cplx G)$.  Since $(H\circ \Phi)(j_{s\ast}\cplx G)\simeq j_{s\ast} (H_s\circ \Phi_s)(\cplx G)$ by \eqref{e:directimage} and $j_s$ is a closed immersion, the unit morphism $\cplx G\to  (H_s\circ \Phi_s)(\cplx G)$ is an
isomorphism; this proves that $\Phi_s$ is fully faithful.

Now assume that $\Phi_s$ is fully faithful for any closed point
$s\in S$. Let us see that the unit
morphism $\eta\colon\Id\to H\circ\Phi$ is an isomorphism. For each
$\cplx F\in \cdbc{X}$we have an exact triangle
$$
\cplx F\xrightarrow{\eta(\cplx F)} (H\circ \Phi)(\cplx F)\to \cono(\eta (\cplx F))\to \cplx F[1]\,.
$$
Then, by \eqref{e:basechange}, for every closed point $s$ in $S$ we have an exact triangle
$$
\bL j_s^\ast \cplx F\to (H_s\circ \Phi_s)(\bL j_s^\ast \cplx F)\to \bL
j_s^\ast \cono(\eta (\cplx F))\to \bL j_s^\ast \cplx F[1]\,.
$$
so that  $\bL j_s^\ast [\cono(\eta (\cplx F))]\simeq \cono(\eta_s(\bL
j_s^\ast \cplx F))\simeq 0$ because $\eta_s\colon \Id \to H_s\circ \Phi_s$. We finish  by Lemma \ref{l:fibres}.

A similar argument gives the statement about equivalence.
\end{proof}

As a corollary of the
previous theorem and Theorem \ref{1:ffcritGo}, we obtain the
following result.

\begin{thm} Let $p\colon X\to S$ and $q\colon Y\to S$ be locally projective Gorenstein morphisms with
integral
%irreducible
 fibres. Let $\cplx{K}$ be an object in
$\cdbc{X\times_SY}$ of finite projective dimension over each
factor. The kernel $\cplx{K}$ is relatively strongly simple over
$X$ (resp.~over $X$ and $Y$) if and only if the functor
$\Phi=\fmf{\cplx{K}}{X}{Y}\colon \cdbc{X}\to \cdbc{Y}$ is fully
faithful (resp.~an equivalence).
\end{thm}

\subsection{Application to Weierstrass elliptic fibrations}

In this subsection we work over an algebraically closed field of characteristic zero.

Let $p\colon X\to S$ be a relatively integral elliptic fibration,
that is, a proper flat morphism whose fibres are integral
Gorenstein curves with arithmetic genus 1. Generic fibres of $p$
are smooth elliptic curves, and the degenerated fibers are
rational curves with one node or one cusp. If $\hat{p}\colon
\hat{X}\to S$ denotes the dual elliptic fibration, defined as the
relative moduli space of torsion free rank 1 sheaves of relative
degree 0, it is known that for every closed point $s\in S$ there
is an isomorphism $\hat{X}_s\simeq X_s$ between the fibers of both
fibrations. If we assume that the original fibration $p\colon X\to
S$ has a section $\sigma\colon S\hookrightarrow X$ taking values
in the smooth locus of $p$, then $p$ and $\hat{p}$ are globally
isomorphic. Let us identify from now on $X$ and $\hat{X}$ and
consider the commutative diagram

$$\xymatrix{ & X\times_S X \ar[ld]_{\pi_1}\ar[rd]^{\pi_2} \ar[dd]^{\rho}&  \\
X \ar[rd]^{p} & & X\ar[ld]_{p} \\
& S & }$$

The relative Poincar\'e sheaf is
$$\cP=\cI_\Delta\otimes \pi_1^\ast\cO_X(H)\otimes \pi_2^\ast\cO_X(H)\otimes
\rho^\ast\omega^{-1}\, ,$$
where $H=\sigma(S)$ is the image of
the section and $\omega=R^1p_\ast\cO_X\simeq (p_\ast \omega_{X/S})^{-1}$.

Relatively integral elliptic fibrations have a Weierstrass form \cite[Lemma II.4.3]{Mir89}: The line
bundle $\cO_X(3H)$ is relatively very ample and if
$\cE=p_\ast \cO_X(3H)\simeq\cO_S\oplus\omega^{\otimes
2}\oplus\omega^{\otimes3}$ and $\bar p\colon \bbP(\cE^\ast)=\operatorname{Proj}(S^\bullet(\cE))\to S$
is the associated projective bundle, there is a closed immersion $j\colon X\hookrightarrow \bbP(\cE^\ast)$
of $S$-schemes such that $j^\ast\cO_{\bbP(\cE^\ast)}(1)=\cO_X(3H)$. In particular, $p$ is a projective morphism.

\begin{lem} The relative Poincar{\'e} sheaf $\cP$ is of finite projective dimension and relatively strongly
simple over both factors.
\end{lem}

\begin{proof} By the symmetry of the expression of $\cP$ it is enough to prove that $\cP$ is of finite
projective dimension and  strongly simple over the first factor.
For the first claim, it is enough to prove that $\cI_\Delta$ has
finite projective dimension. Let us consider the exact sequence
$$0\to \cI_\Delta\to \cO_{X\times_S X}\to \delta_\ast \cO_X\to 0$$
where $\delta\colon X\hookrightarrow X\times_SX$ is the diagonal
morphism. It suffices to see that $\delta_\ast \cO_X$ has finite
projective dimension. We have to prove that for any $\cplx{N}\in
\dbc{X}$, the complex $\bR\dSHom{\cO_{X\times_S X}} (\delta_\ast
\cO_X, \pi_1^!\cplx{N})$ is bounded. This is a complex supported
at the diagonal, so that  it suffices to see that
$\bR\pi_{1\ast}\bR\dSHom{\cO_{X\times_S X}} (\delta_\ast \cO_X,
\pi_1^!\cplx{N})$ is bounded. This follows from the following formulas.
$$
 \bR\pi_{1\ast}\bR\dSHom{\cO_{X\times_S X}}(\delta_\ast \cO_X,
\pi_1^!\cplx{N})
 \simeq \bR\dSHom{\cO_{X}} ( \cO_X,
\cplx{N})\simeq\cplx{N}\,.
$$

Let us prove that $\cP$ is strongly simple over the first factor.  Fix a
closed point $s\in S$ and consider two different points $x_1$ and $x_2$ in
the fiber $X_s$.  If both are non-singular, then
$$
\Hom^i_{\catD(X_s)}(\fmf{\cP_s}{X_s}{X_s}(\cO_{x_1}),\fmf{\cP_s}{X_s}{X_s}(\cO_{x_2}))\simeq
H^i(X_s,\cO_{X_s}(x_1-x_2))=0\quad \text{for every $i$}
$$
because $\cO_{X_s}(x_2-x_1)$ is a non-trivial line bundle of degree zero. Assume that $x_2$ is
singular and $x_1$ is not, the other case being similar. Let $Z_{x_2}$ be  a l.c.i
zero cycle supported on $x_2$. We have
$$
\Hom^i_{\catD(X_s)}(\fmf{\cP_s}{X_s}{X_s}(\cO_{x_1}),\fmf{\cP_s}{X_s}{X_s}(\cO_{Z_{x_2}}))=
H^i(X_s, \cJ_{Z_{x_2}}\otimes \cO_{X_s}(x_1))
$$
where $\cJ_{Z_{x_2}}$ denotes the direct image by the finite morphism $Z_{x_2}\times X_s\to X_s$ of
the ideal sheaf of the graph $Z_{x_2}\hookrightarrow Z_{x_2}\times X_s$ of the immersion
$Z_{x_2}\hookrightarrow X_s$.

Let us consider the exact sequences of $\cO_{X_s}$-modules
\begin{align*}
0\to \cJ_{Z_{x_2}} \to \cO_{Z_{x_2}}&\otimes_k \cO_{X_s}\to \cO_{Z_{x_2}} \to 0 \\
0\to \cJ_{Z_{x_2}}(x_1) \to \cO_{Z_{x_2}}&\otimes_k \cO_{X_s}(x_1)\to \cO_{Z_{x_2}} \to 0
\end{align*}
Since $H^0(X_s,\cO_{X_s})\simeq k$  the morphism $\cO_{Z_{x_2}}\otimes_k H^0(X_s, \cO_{X_s})\to
\cO_{Z_{x_2}}$ of global sections induced by the first sequence is an isomorphism. Moreover,
$H^0(X_s,\cO_{X_s})\simeq H^0(X_s,\cO_{X_s}(x_1))$ and then  we also have an isomorphism of global
sections $\cO_{Z_{x_2}}\otimes_k H^0(X_s, \cO_{X_s})\iso\cO_{Z_{x_2}}\otimes_k H^0(X_s, \cO_{X_s}(x_1))$.
Thus, $\cO_{Z_{x_2}}\otimes_k H^0(X_s, \cO_{X_s}(x_1))\iso \cO_{Z_{x_2}}$ so that
$H^i(X_s,\cJ_{Z_{x_2}}(x_1))=0$ for $i=0,1$.

Finally, since $\Hom^0_{\catD(X_s)}(\cP_x, \cP_x)= k$ for every point
$x\in X_s$, we conclude that $\cP_s$ is strongly simple over
$X_s$.
\end{proof}
Now by Corollary \ref{c:ffcritGo} we have
\begin{prop} \label{p:poincequiv}
The relative integral functor
$$\fmf{\cP}{X}{X}\colon \cdbc{X}\to \cdbc{X}$$  defined by the Poincar{\'e} sheaf is an equivalence of categories.
\end{prop}
Notice that the proof of this result does not use spanning classes.

%\backmatter

\end{document}